%% file: lsat.tex
\documentclass[12pt]{article}
\usepackage{fullpage}
\usepackage{booktabs}
\usepackage{amsmath}
\usepackage{siunitx}
\usepackage{hyperref}
\usepackage{graphicx}
\usepackage{subcaption}
\input{defs.tex}

\begin{document}

\title{Least Squares Auto-Tuning}
\author{
Shane Barratt \\
sbarratt@stanford.edu
\and
Stephen Boyd \\
boyd@stanford.edu
}

\date{}

\maketitle

\begin{abstract}
Least squares is by far the simplest and most commonly applied
computational method in many fields.
In almost all applications, the least squares objective is rarely the true objective.
We account for this discrepancy by parametrizing the least squares problem
and automatically adjusting these parameters using an optimization algorithm.
We apply our method, which we call least squares auto-tuning, to data fitting.
\end{abstract}

\section{Introduction}

Since its introduction over 200 years ago by Legendre and Gauss,
the method of least squares~\cite{legendre1805nouvelles,gauss1809theoria}
has been one of the most widely employed computational techniques in
many fields, including machine learning and statistics, signal processing,
control, robotics, and finance \cite{boyd2018introduction}.
Its wide application primarily comes from the fact that it has a simple
analytical solution, it is easy to understand, and very efficient and
stable algorithms for computing its solution have been developed
\cite{lawson1995solving,golub2012matrix}.

In essentially all applications, the least squares objective is not
the true objective; rather it is a surrogate for the real goal.
For example, in least squares data fitting, the objective is not to solve 
a least squares problem involving the training data set, 
but rather to find a model or predictor that generalizes, \ie, achieves small
error on new unseen data.
In control, the least squares objective is only a surrogate for keeping the
state near some target or desired value, while keeping
the control or actuator input small.

To account for the discrepancy between the least squares objective and the
true objective, it is common practice to modify (or tune) the least squares
problem that is solved to obtain a good solution in terms of the true objective.
Typical tricks here include modifying the data, adding additional 
(regularization) terms to the cost function, or
varying hyper-parameters or weights in the least squares problem to be solved.

The art of using least squares in applications is generally in
how to carry out these modifications or choose these additional terms, and how
to choose the hyper-parameters. The choice of hyper-parameters is often done in
an ad hoc way, by varying them, solving the least squares problem, and then
evaluating the result using the true objective or objectives.
In data fitting, for example, regularization scaled by a hyper-parameter
is added to the least squares problem, which is solved for many values of the
hyper-parameter to obtain a family of data models;
among these, the one that gives the best predictions on a
test set of data is the one that is ultimately used.
We refer to this general design approach, of modifying
the least squares problem to be solved,
varying some hyper-parameters, and evaluating the result using the
true objective, as \emph{least squares tuning}.
It is very widely used, and can be extremely effective in practice.

Our focus in this paper is on automating the process of least squares
tuning, for a variety of data fitting applications.
We parametrize the least squares problem to be solved by hyper-parameters,
and then automatically adjust these hyper-parameters using a
gradient-based optimization algorithm, to obtain the best
(or at least better) true performance.
This lets us automatically search the hyper-parameter design space,
which can lead us to better designs than could be found manually, or help us
find good values of the hyper-parameters more quickly than if the adjustments
were done manually.
We refer to the method as \emph{least squares auto-tuning}.

One of our main contributions in this paper is the observation that
least squares auto-tuning is very effective for a wide variety of data
fitting problems that are usually handled using more complex and
advanced methods, such as non-quadratic loss functions or regularizers
in regression, or special loss functions for classification problems.
In addition, it can simultaneously adjust hyper-parameters in the
feature generation chain. Through several examples, we show that ordinary
least squares, used for over 200 years, coupled with automated hyper-parameter
tuning, can be very effective as a method for data fitting.

The method we describe for least squares auto-tuning is easy to
understand and just as easy to implement.
Moreover, it is an exercise in calculus
to find the derivative of the least squares solution, and an exercise in
numerical linear algebra to compute it efficiently.
We describe an implementation that utilizes new and powerful software frameworks
that were originally designed to optimize the parameters in deep
neural networks, making it very efficient on modern hardware and
allowing it to scale to (extremely)
large least squares tuning problems.

\paragraph{Our contributions.}
We claim three main contributions.
The first contribution is the observation that the least squares
solution map can be efficiently differentiated, including when
the problem data is sparse; we mirror our description with an
open-source implementation for both of these cases.
The second contribution is the method of least squares auto-tuning,
which can automatically tune hyper-parameters in least squares problems.
The final contribution is our unique application of least squares auto-tuning
to data fitting.

\section{Background and related work}

Our work mainly falls at the intersection of two fields: automatic differentiation
and hyper-parameter optimization.
In this section we review related work.

\paragraph{Automatic differentiation.}
The general idea of automatic differentiation (AD) is to automatically compute
the derivatives of a function given a program that evaluates the function
\cite{wengert1964simple,speelpenning1980compiling}. In general, the cost of
computing the derivative or gradient of a function can be made about the same
(usually within a factor of 5) as computing the function
\cite{baur1983complexity,griewank2008evaluating}. This means that an
optimization algorithm can obtain derivatives of the function it is optimizing
as fast as computing the function itself, and explains the proliferation
of gradient-based minimization methods
\cite{baydin2018automatic,bottou2018optimization}.
There are many popular implementations of AD, and they generally fall
into two categories. The first category is are trace-based AD systems, which
trace computations at runtime, as they are executed;
popular ones include PyTorch \cite{paszke2017automatic},
Tensorflow eager \cite{agrawal2019tensorflow},
and autograd \cite{maclaurin2015autograd}.
The second category are based on source transformation, which transform
the (native) source code that implements the function into source code
that implements the derivative operation. Popular implementations here
include Tensorflow \cite{abadi2016tensorflow}, Tangent \cite{van2018tangent},
and (more recently) Zygote \cite{innes2018don}.

\paragraph{Argmin differentiation.}
Given an optimization problem parametrized by some parameters, the
solution map is a set-valued map from those parameters to a set of solutions.
If the solution map is differentiable (and in turn unique), then we can differentiate
the solution map~\cite{dontchev2009implicit}.
For convex optimization problems that satisfy strong duality, the
solution map is given by the set of solutions to the KKT conditions,
which can in some cases be differentiated using the implicit function theorem
\cite{barratt2018differentiability}. This idea has been applied to convex
quadratic programs \cite{amos2017optnet}, stochastic optimization
\cite{donti2017task}, games \cite{ling2018game,ling2019large},
physical systems \cite{de2018end}, control \cite{amos2018differentiable},
and structured inference \cite{belanger2016structured,belanger2017end}.
In machine learning, these
techniques were originally applied to neural networks
\cite{larsen1998adaptive,eigenmann1999gradient} and ridge regression
\cite{bengio2000gradient}, and more recently to lasso \cite{mairal2012task},
support vector machines~\cite{chapelle2002choosing},
and log-linear models \cite{keerthi2007efficient,foo2008efficient}
A notable AD implementation of these methods is the PyTorch implementation
\verb|qpth|,
which can compute derivatives of the solution map of quadratic programs
\cite{qpth}.
\nocite{luketina2016scalable}
\nocite{pedregosa2016hyperparameter}

\paragraph{Unrolled optimization.}
Another approach to argmin differentiation is unrolled optimization.
In unrolled optimization, one fixes the number of iterations in an
iterative minimization method, and differentiates the steps taken by the
method itself \cite{domke2012generic,baydin2014automatic}.
The idea of unrolled optimization was originally applied to optimizing
hyper-parameters in deep neural
networks, and has been extended in several ways to adjust learning rates,
regularization parameters
\cite{maclaurin2015gradient,fu2016distilling,lorraine2018stochastic},
and even to learn weights on individual data points \cite{ren2018learning}.
It is still unclear whether argmin differentiation should be performed
via implicit differentiation or unrolled optimization.
However, when the optimization problem is nonconvex, differentiation
by unrolled optimization seems to be the only practical one.

\paragraph{Hyper-parameter optimization.}
The idea of adjusting hyper-parameters to obtain better true performance in
the context of data
fitting is hardly new, and routinely employed in settings more sophisticated
than least squares.
For example, in data fitting, it is standard practice to
vary one or more hyper-parameters to generate a set of models,
and choose the model that attains the best true objective,
which is usually error on an unseen test set.
The most commonly employed methods here include grid search, random search
\cite{bergstra2012random}, Bayesian optimization 
\cite{movckus1975bayesian,rasmussen2004gaussian,snoek2012practical},
and covariance matrix adaptation \cite{hansen1996adapting}.

\section{Least squares auto-tuning}
\label{sec:lsat}

In this section we describe the idea of least squares tuning,
our method for least squares auto-tuning, and
our implementation.

\subsection{Least squares problem}

The \emph{matrix least squares problem} that depends on a
\emph{hyper-parameter vector} $\omega\in \Omega \subseteq \reals^p$
has the form
\BEQ
\begin{array}{ll}
\mbox{minimize} & \|A(\omega)\theta-B(\omega)\|_F^2,
\end{array}
\label{eq:lstsq}
\EEQ
where the variable is $\theta\in\reals^{n \times m}$, 
the \emph{least squares optimization variable}
or \emph{parameter matrix},
and $A: \Omega \rightarrow \reals^{k \times n}$ and 
$B: \Omega \rightarrow \reals^{k \times m}$ map the hyper-parameter vector
to the least squares \emph{problem data}.
The norm $\| \cdot \|_F$ denotes the Frobenius norm, \ie, the squareroot of
the sum of squares of the entries of a matrix.
We assume throughout this paper that $A(\omega)$ 
has linearly independent columns, which implies that it is tall, \ie,
$k \geq n$.
Under these assumptions, the \emph{least squares solution} is unique, given by
\BEQ
\theta^{\mathrm{ls}}(\omega) = 
A(\omega)^\dagger B(\omega) =
(A(\omega)^TA(\omega))^{-1}A(\omega)^T B(\omega),
\label{eq:lstsqsolution}
\EEQ
where $A(\omega)^\dagger$ denotes the (Moore-Penrose) pseudo-inverse.
Solving a least squares problem for a given hyper-parameter vector 
corresponds to computing $\theta^\mathrm{ls}(\omega)$.
We will think of the least squares solution $\theta^\mathrm{ls}$ as a
function mapping the hyper-parameter $\omega \in \Omega$ to a parameter 
$\theta^\mathrm{ls}(\omega) \in \reals^{n \times m}$.

\paragraph{Multi-objective least squares.}
In many applications we have multiple least squares objectives.  
These are typically scalarized
by forming a positive weighted sum, which leads to
\BEQ
\begin{array}{ll}
\mbox{minimize} & \lambda_1 \|A_1(\omega)\theta-B_1(\omega)\|_F^2 + \cdots +
				  \lambda_r \|A_r(\omega)\theta-B_r(\omega)\|_F^2,
\end{array}
\label{eq:lstsq-multi}
\EEQ
where $\lambda_1, \ldots, \lambda_r$ are the positive objective weights.
This problem is readily expressed as the standard least squares
problem~(\ref{eq:lstsq}) by stacking the objectives, with
\BEQ
A(\omega) =
\begin{bmatrix}
\sqrt{\lambda_1} A_1(\omega) \\ \vdots \\ \sqrt{\lambda_r} A_r(\omega)
\end{bmatrix},
\quad B(\omega) =
\begin{bmatrix}
\sqrt{\lambda_1} B_1(\omega) \\ \vdots \\ \sqrt{\lambda_r} B_r(\omega)
\end{bmatrix}.
\label{eq:transformation}
\EEQ
We will often write least squares problems in the form \eqref{eq:lstsq-multi}, 
and assume that the reader understands that the problem data can
easily be transformed into \eqref{eq:transformation}.
The objective weights $\lambda_1, \ldots, \lambda_r$ can also be considered 
hyper-parameters themselves, or to depend on hyper-parameters;
to keep the notation light we do not show this dependence.

\paragraph{Solving the least squares problem.}
For a given value of $\omega$, there are many ways to solve the
least squares problem \eqref{eq:lstsq},
including dense or sparse QR or other factorizations
\cite{golub1965numerical,bjorck1980direct}, 
iterative methods such as CG or LSQR
\cite{hestenes1952methods,paige1982lsqr}, and
many others \cite{lawson1995solving,golub2012matrix}.
Very efficient libraries for computing the least squares solution that
target multiple CPUs or one
or more GPUs have also been developed
\cite{dongarra1990algorithm,anderson1999lapack,walt2011numpy,sanders2010cuda}.
We note that the problem is separable across the columns of $\theta$, \ie,
the problem splits into $m$ independent least squares problems with
vector variables and a common coefficient matrix.

We give a few more details here for two of these methods.
First we consider the case where $A(\omega)$ and
$B(\omega)$ are stored and manipulated as dense matrices.
One attractive option (for GPU implementation) is to form the
Gram matrix $G=A^TA$, along with $H=A^TB$.
This requires around (order) $kn^2$ and $knm$ flops, respectively, but
these matrix-matrix multiplies are BLAS level 3 operations,
which can be carried out very efficiently.
To compute $\theta^\mathrm{ls}$, we can use a Cholesky factorization of
$G$, $G=LL^T$, which costs
order $n^3$ flops, solve the triangular equation $LY=H$, which costs
order $n^2m$ flops, and then solve the triangular equation
$L^T \theta^\mathrm{ls} = Y$, which costs order $n^2m$ flops.
Overall, the complexity of solving a dense least squares problem
is order $kn(n+m)$.

The other case for which we give more detail is when
$A(\omega)$ is represented as an abstract linear operator, and not as a matrix.
This is a natural representation when $A(\omega)$ is large and sparse,
or represented as a product of small (or sparse) matrices.
That is, we can evaluate $A(\omega)u$ for any $u\in \reals^n$, and
$A(\omega)^T v$ for any $v\in \reals^k$.
(This is the so-called \emph{matrix-free} representation.)
We can use CG or LSQR to solve the least squares problem, in
parallel for each column of $\theta$. The complexity of CG or LSQR
depends on the problem, and can vary considerably 
based on the data, size, sparsity, choice 
of pre-conditioner, and required accuracy
\cite{hestenes1952methods,paige1982lsqr}.

\subsection{Least squares tuning problem}
In a \emph{least squares tuning problem}, our goal is to choose
the hyper-parameters to achieve some goal. We formalize this as the problem
\BEQ
\begin{array}{ll}
\mbox{minimize} & F(\omega)
=
\psi(\theta^\mathrm{ls}(\omega)) + r(\omega),
\end{array}
\label{eq:lsep}
\EEQ
with variable $\omega \in \Omega$ and objective
$F:\Omega\rightarrow \reals \cup \{+\infty\}$, where
$\psi: \reals^{n \times m} \rightarrow \reals$ is the 
\emph{true objective function}, and
$r: \Omega \rightarrow \reals \cup \{+\infty\}$ is the 
\emph{hyper-parameter regularization function}.
We use infinite values of $r$ (and therefore $F$) to encode constraints 
on the hyper-parameter $\omega$, and will assume that
$r(\omega)$ is defined as $\infty$
for $\omega \not\in \Omega$.
A least squares tuning problem is specified by the functions 
$A$, $B$, $\psi$, and $r$.
We will make some additional assumptions about these functions below.

The hyper-parameter regularization function $r$ can itself contain
a few parameters that can be varied, which of course affects
the hyper-parameters chosen in the least squares auto-tuning
problem~(\ref{eq:lsep}), which in turn affects the parameters
selected by least squares.  We refer to parameters that may appear 
in the hyper-parameter regularization function as \emph{hyper-hyper-parameters}.

The least squares tuning problem~(\ref{eq:lsep}) can be formulated in
several alternative ways, for example as the constrained problem
with variables $\theta \in \reals^{n \times m}$ 
and $\omega \in \Omega$
\BEQ
\begin{array}{ll}
\mbox{minimize} & \psi(\theta) + r(\omega) \\
\mbox{subject to} & A(\omega)^TA(\omega) \theta = A(\omega)^TB(\omega).
\end{array}
\label{eq:probconstrained}
\EEQ
In this formulation, $\theta$ and $\omega$ are independent variables, coupled
by the constraint, which is the optimality condition for the least squares
problem~(\ref{eq:lstsq}).
Eliminating the constraint in this problem yields our formulation~(\ref{eq:lsep}).

\paragraph{Solving the least squares tuning problem.}
The least squares tuning problem is in general nonconvex, and difficult
or impossible as a practical matter to solve exactly
\cite{polyak1987introduction,boyd2004convex}.
(One important exception
is when $\omega$ is a scalar and $\Omega$ is an interval,
in which case we can simply evaluate 
$F(\omega)$ on a grid of values over $\Omega$.)
This means that we will need to resort to a local optimization
or heuristic method in order to (approximately) solve it.

We will assume that $A$ and $B$ are differentiable in $\omega$,
which implies that $\theta^\mathrm{ls}$ is differentiable in $\omega$,
since the mapping from $\omega$ to $\theta^\mathrm{ls}$ is differentiable.
We will also assume that $\psi$ is differentiable, which implies 
that the true objective $\psi(\theta^\mathrm{ls}(\omega))$
is differentiable in the hyper-parameters $\omega$. 
This means that the first term (the true objective) in the least
squares tuning problem \eqref{eq:lsep} is
differentiable, while the second one (the hyper-parameter regularizer)
need not be. There are many methods that can be used
to (approximately) solve such a composite problem
\cite{douglas1956numerical,lions1979splitting,shor1985minimization,
boyd2011distributed,nesterov2013introductory,parikh2014proximal}.

For completeness, we describe one of the simplest methods, 
the proximal gradient method (which stems from the proximal
point method \cite{martinet1970breve};
for a modern reference see \cite{nesterov2013gradient}), given by the iteration
\[
\omega^{k+1} =
\mathrm{\bf prox}_{t^k r}
(\omega^k - t^k \nabla_\omega \psi(\theta^\mathrm{ls}(\omega^k))),
\]
where $k$ denotes the iteration number, $t^k>0$ is a step size, and
the proximal operator $\mathrm{\bf prox}_{tr}: \reals^p \rightarrow \Omega$
is given by
\[
\mathrm{\bf prox}_{tr}(\nu) = \argmin_{\omega \in \Omega}
\left(t r(\omega) + (1/2)\|\omega-\nu\|_2^2\right).
\]
We assume here that the argmin exists; when it is not unique, we choose
any minimizer.
In order to use the proximal gradient method, we
need the proximal operator of $tr$ to be relatively easy to evaluate.

The proximal gradient method reduces to the ordinary gradient method when
$r=0$ and $\Omega = \reals^p$.
Another special case is when $\Omega \subset \reals^p$, and
$r(\omega)=0$ for $\omega\in \Omega$.
In this case $r$ is the indicator function of the set $\Omega$,
the proximal operator of $tr$ is Euclidean projection onto $\Omega$, and the 
proximal gradient method coincides with the projected gradient method.

\paragraph{Choosing the step size.}

There are many ways to choose the step size.
We adopt the following simple adaptive scheme,
borrowed from \cite{lall2017ee104}.
The method begins with an initial step size $t^1$.
If the function value decreases or stays the same from iteration
$k$ to $k+1$, or $F(\omega^{k+1}) \leq F(\omega^k)$,
then we increase the step size a bit 
and accept the update.
If, on the other hand, the function value increases from iteration
$k$ to $k+1$, or
$F(\omega^{k+1}) > F(\omega^k)$, 
we decrease the step size substantially and reject the update,
\ie, $\omega^{k+1} = \omega^k$.
A simple rule for increasing the step size is $t^{k+1}=(1.2)t^k$,
and a simple rule for decreasing it is $t^{k+1} = (1/2) t^k$.

We note that more sophisticated step size selection methods exist
(see, \eg, the line search methods for the Goldstein or Armijo conditions
\cite{nocedal2006numerical}).

\paragraph{Stopping criterion.}

By default, we run the method for a fixed maximum number of iterations.
If $F(\omega^{k+1}) \leq F(\omega^k)$, then a reasonable stopping
criterion at iteration $k+1$ is
\BEQ
\|(\omega^k - \omega^{k+1})/t^k + (g^{k+1} - g^k) \|_2 \leq \epsilon,
\label{eq:stopping-criterion}
\EEQ
where $g^k = \nabla_{\omega^k} \psi(\theta^\mathrm{ls}(\omega^k))$,
for some small tolerance $\epsilon > 0$.
(For more justification of this stopping criterion,
see Appendix \ref{sec:deriv_stopping_criterion}.)
When $r=0$ and $\Omega = \reals^p$
(\ie, the proximal gradient method coincides with the ordinary gradient method),
this stopping criterion reduces to
\[
\|g^{k+1}\|_2 \leq \epsilon,
\]
which is the standard stopping criterion in the ordinary gradient method.
The full algorithm for least squares auto-tuning via the proximal gradient method
is summarized in Algorithm \ref{alg:lstsq}.

\begin{algdesc}
\label{alg:lstsq}
\emph{Least squares auto-tuning via proximal gradient.}
\begin{tabbing}
    {\bf given} initial hyper-parameter vector $\omega^1 \in \Omega$,
    initial step size $t^1$, number of iterations $n_\mathrm{iter}$,\\
    \qquad \=\ tolerance $\epsilon$.\\
    {\bf for} $k=1,\ldots,n_\mathrm{iter}$\\
	    \qquad \=\ 1.\ \emph{Solve the least squares problem}. 
	    $\theta^\mathrm{ls}(\omega^k) = (A^T(\omega^k)A(\omega^k))^{-1}A^T(\omega^k)B(\omega^k)$.\\
	    \qquad \=\ 2.\ \emph{Compute the gradient}.
		$g^k = \nabla_\omega \psi(\theta^\mathrm{ls}(\omega^k))$.\\
	    \qquad \=\ 3.\ \emph{Compute the gradient step}. $\omega^{k+1/2} = \omega^k - t^k g^k$.\\
	    \qquad \=\ 4.\ \emph{Compute the proximal operator}. 
		$\omega^\mathrm{tent} = \mathrm{\bf prox}_{t^k r} (\omega^{k+1/2})$.\\
	    \qquad \=\ 5. {\bf if} $F(\omega^\mathrm{tent}) \leq F(\omega^k)$,\\
	    \qquad \qquad \=\ \emph{Increase step size and accept update.}
	    $t^{k+1} = (1.2)t^k; \quad \omega^{k+1} = \omega^\mathrm{tent}$.\\
	    \qquad \qquad \=\ \emph{Stopping criterion.} {\bf quit} if
	    $\|(\omega^k - \omega^{k+1})/t^k + (g^{k+1} - g^k) \|_2 \leq \epsilon$.\\
	    \qquad \=\ 6. {\bf else} \emph{Decrease step size and reject update}.
	    $t^{k+1} = (1/2)t^k; \quad \omega^{k+1} = \omega^k$.\\
    {\bf end for}
\end{tabbing}
\end{algdesc}

We emphasize that many other methods can be used to (approximately)
solve the least squares tuning problem \eqref{eq:lsep};
we have described the proximal gradient method here only for completeness.

\subsection{Computing the gradient}
\label{sec:computing_the_gradient}

\emph{Note.} In principle, one could calculate the gradient by directly
differentiating the linear algebra routines used to solve the least
squares problem \cite{smith1995differentiation}.
We, however, work out formulas for computing the gradient analytically, that
work in the case of sparse $A$ and allow for a more efficient implementation.

To compute $g=\nabla_\omega \psi(\theta^\mathrm{ls}(\omega)) \in \reals^p$
we can make use of the chain rule
for the composition $f=\psi \circ \theta^\mathrm{ls}$.
We assume that $\theta = \theta^\mathrm{ls}(\omega)$ has been computed.
We first compute $\nabla_{\theta}\psi(\theta) \in \reals^{n \times m}$,
and then form
\[
C=(A^TA)^{-1}\nabla_\theta \psi(\theta) \in \reals^{n \times m}.
\]
(Here we have dropped the dependence on $\omega$, \ie, $A=A(\omega)$.)
If $A$ is stored as a dense matrix, we observe that $G=A^TA$ 
and its factorization have already been computed (to evaluate $\theta$), so
this step involves a back-solve.
If $A$ is represented as an abstract operator,
we can evaluate each column of $C$ (in parallel) using an iterative method.

It can be shown (see Appendix \ref{sec:deriv_grads})
that the gradients of $\psi$ with respect to $A$ and $B$ are given by
\BEQ
\nabla_A \psi = (B-A\theta)C^T-AC\theta^T \in \reals^{k \times n}, \qquad
\nabla_B \psi = AC \in \reals^{k \times m}.
\label{eq:grads}
\EEQ
(Again, the dependence on $\omega$ has been dropped.)

In the case of $A$ dense, we can explicitly form $\nabla_A \psi$ and
$\nabla_B \psi$, since they are
the same size as $A$ and $B$, which we already have stored.
The overall complexity of computing
$\nabla_A \psi$ and $\nabla_B \psi$ is $kn(n+m)$ in the dense case,
which is the same cost
as solving the least squares problem.

In the case of $A$ sparse, we can explicitly form the matrix $\nabla_B \psi$,
but we can not form
$\nabla_A \psi$, since it is the size of $A$, 
which by assumption is too large to store.
Instead, we assume that $\omega$ only affects $A$ at a subset of its entries,
$\Gamma$, \ie, $A_{ij}(\omega)=0$
for all $i,j \not\in \Gamma$,
and for all $\omega \in \Omega$.
By doing this, we have restricted $\nabla_A \psi$ to have the
same sparsity pattern as $A$, meaning we only need to compute
$(\nabla_A \psi)_{ij}$ for $i,j \in \Gamma$.
That is, we compute
\[
\nabla_A \psi = \begin{cases}
(b_i-\theta a_i) c_j^T - a_i^T (C\theta^T)_j & i,j \in \Gamma \\
0 & \text{otherwise},
\end{cases}
\]
where $b_i$ is the $i$th row of $B$, $a_i$ is the $i$th row of $A$,
$c_j$ is the $j$th column of $C$, and $(C\theta^T)_j$ is the $j$th column
of $C\theta^T$.
(This computation can be done in parallel.)

The next step is to compute $g=\nabla_\omega \psi$ given $\nabla_A \psi$ and $\nabla_B \psi$.
We first describe how $g$ can be computed in the case of dense $A$, and then in the case of sparse $A$.

\paragraph{Dense $A$.}

We first evaluate
$\nabla_{\omega} A_{ij} \in \reals^p$ and $\nabla_{\omega} B_{ij}\in \reals^p$, 
the gradients of the problem data entries with respect to $\omega$.
If these gradients are all dense, we need to store $k(n+m)$ vectors in $\reals^p$;
but generally, they are quite sparse.
(We explain how to take advantage of the sparsity in \S \ref{sec:implementation}.)
Finally, we have
\[
g = \sum_{i,j}
(\nabla_A \psi)_{ij}
(\nabla_{\omega} A)_{ij}
+ \sum_{i,j}
(\nabla_B \psi)_{ij}
(\nabla_{\omega} B)_{ij}.
\]
Assuming these are all dense, this requires order $knp+kmp = kp(n+m)$ flops.
The overall complexity of
evaluating the gradient $g$ is order $k(n+m)(n+p)$.

\paragraph{Sparse $A$.}

When $A$ is large and sparse, we only need to compute the inner product
at the entries of $A$ that are affected by $\omega$, \ie, we compute
\[
g = \sum_{i,j \in \Gamma}
(\nabla_A \psi)_{ij}
(\nabla_{\omega} A)_{ij}
+ \sum_{i,j}
(\nabla_B \psi)_{ij}
(\nabla_{\omega} B)_{ij}.
\]
If $|\Gamma| \ll kn$, then this can be much faster than treating
$A$ as dense.

\subsection{Implementation}
\label{sec:implementation}

The equations in \S \ref{sec:computing_the_gradient} for computing $g$
do not directly lend themselves to an implementation.
For example, we need to compute $\nabla_\theta \psi$, $\nabla_\omega A_{ij}$ and $\nabla_\omega B_{ij}$,
which depend on the form of $\psi$, $A$, and $B$.
Also, we would like to take advantage of the (potential) sparsity of these gradients.
We can use libraries for automatic differentiation, \eg, PyTorch \cite{paszke2017automatic}
and Tensorflow \cite{abadi2016tensorflow},
to automatically (and efficiently) compute $g$ given $\psi$, $A$, and $B$.

These libraries generally work by representing functions as differentiable computation graphs,
allowing one to evaluate the function as well as its gradient.
In our case, we represent $\psi$ as a function of $\omega$, defined by a differentiable computation graph,
and use these libraries to automatically compute $g=\nabla_\omega \psi$.
In order to use these libraries to compute $g$, we need to implement an operation that solves
the least squares problem \eqref{eq:lstsqsolution} and computes its gradients \eqref{eq:grads}.
We have implemented an operation \verb+lstsq(A,B)+ that does exactly this, in both PyTorch and Tensorflow, in both the dense and sparse case.
(The code can be found in the Appendix.)

There are several advantages of using these libraries.
First, they automatically exploit parallelism and gradient sparsity.
Second, they utilize BLAS level 3 operations, which are very efficient on modern hardware.
Third, they make it easy to represent the functions $\psi$, $A$, and $B$, since
they can be represented as compositions of (the many) pre-defined operations
in these libraries.

We also provide a generic PyTorch implementation of the adaptive
proximal gradient algorithm that we used in this paper.

\paragraph{GPU timings.}

Table \ref{tab:gpu} gives timings of computing
$\psi(\theta^\mathrm{ls}(\omega))=\Tr(\ones\ones^T\theta^\mathrm{ls}(\omega))$
and its gradient $g$
for a random problem (where $\omega$ simply scales the rows of a fixed $A$ and $B$)
and various problem dimensions, where $A$ is dense.
The timings given are for the PyTorch implementation, on an unloaded
GeForce GTX 1080 Ti Nvidia GPU using 32-bit floating point numbers (floats).
The timings are about ten times longer using 64-bit floating point numbers (doubles).
(For most applications, including data fitting, we only need floats.)
We also give the percentage of time spent on the Cholesky factorization
of the Gram matrix.

\begin{table}
  \caption{GPU timings.}
  \centering
  \begin{tabular}{lllllll}
    \toprule
    $k$ & $n$ & $m$ & $p$ & Compute $\psi$ & Compute $g$ & Cholesky \\
    \midrule
    20000 & 10000 & 10000 & 20000 & \SI{1.32}{\second} & \SI{2.49}{\second} & 15.6 \% \\
    20000 & 10000 & 1 & 20000 & \SI{446}{\milli\second} & \SI{614}{\milli\second} & 16.1 \% \\
    100000 & 1000 & 100 & 100000 & \SI{28}{\milli\second} & \SI{28}{\milli\second} & 10.8 \% \\
    \bottomrule
  \end{tabular}
  \label{tab:gpu}
\end{table}

\subsection{Equality constrained extension}
\label{sec:ext_var}

One can easily extend the ideas described in this paper to the equality-constrained
least squares problem
\[
\begin{array}{ll}
\text{minimize} & \|A(\omega)\theta-B(\omega)\|_F^2 \\
\text{subject to} & C(\omega)\theta = D(\omega),
\end{array}
\]
with variable $\theta\in\reals^{n \times m}$,
where $C:\Omega \to \reals^{d \times n} $ and $D:\Omega \to \reals^{d \times m}$.
The pair of primal and dual variables $(\theta, \nu) \in
\reals^{n \times m} \times \reals^{d \times m}$ are optimal if and only if they
satisfy the KKT conditions \cite[Chapter 16]{boyd2018introduction}
\[
\begin{bmatrix}
A(\omega)^TA(\omega) & C(\omega)^T \\
C(\omega) & 0 
\end{bmatrix}
\begin{bmatrix} \theta \\ \nu \end{bmatrix}
=
\begin{bmatrix}
A(\omega)^T B(\omega) \\ D(\omega)
\end{bmatrix}.
\]
From here on, we let
\[
M(\omega) =
\begin{bmatrix}
A(\omega)^TA(\omega) & C(\omega)^T \\
C(\omega) & 0
\end{bmatrix}.
\]
When $A$ and $C$ are dense, one can factorize $M$ directly,
using an $LDL^T$ factorization \cite{lawson1995solving}.
When $A$ and $C$ are sparse matrices, one can solve this system iteratively
(\ie, without forming the matrix $M$) using, \eg, MINRES \cite{paige1975solution}.

Our true objective function becomes a function of both $\eta=(\theta,\nu)$.
Suppose we have the gradient $\nabla_\eta \psi$.
We first compute
\[
\begin{bmatrix}H_1 \\ H_2 \end{bmatrix} = 
M(\omega)^{-1}
\nabla_\eta \psi
\]
and
\[
F =
M(\omega)^{-1}
\begin{bmatrix}A^TB \\ 0\end{bmatrix}.
\]
Then the gradients of $\psi$ with respect to $A$ and $B$ are given by
\[
\nabla_A \psi = BH_1^T - A(FH_1^T + AH_1 F^T), \quad
\nabla_B \psi = AH_1.
\]
and with respect to $C$ and $D$ are given by
\[
\nabla_C \psi = -\nu H_1^T - H_2 \theta^T, \quad
\nabla_D \psi = H_2.
\]
Computing the gradients of the solution map requires the solution of two linear systems,
and thus has roughly double the complexity of computing the solution itself
(and much less in the dense case when a factorization is cached).
When $A$ and $C$ are sparse, we can compute their gradients at only the nonzero elements,
in a similar fashion to the procedure described in \S\ref{sec:computing_the_gradient}.

\section{Least squares data fitting}
\label{sec:auto_lsd}

In the previous section we described the general idea of least squares auto-tuning.
In this section, and for the remainder of the paper, we apply least squares auto-tuning to
data fitting.

\subsection{Least squares data fitting}
In a data fitting problem, we have \emph{training data} consisting of \emph{inputs}
$u_1,\ldots,u_N\in\mathcal U$ and \emph{outputs} $y_1,\ldots,y_N\in\reals^m$.
In \emph{least squares data fitting}, we fit the parameters of a predictor
\[
\hat y = \phi(u,\omega^\mathrm{feat})^T \theta,
\]
where $\theta\in\reals^{n\times m}$ is the model variable, 
$\omega^\mathrm{feat} \in \Omega^\mathrm{feat} \subseteq \reals^{p^\mathrm{feat}}$ 
are feature engineering hyper-parameters, and 
$\phi: \mathcal U \times \Omega^\mathrm{feat} \to \reals^n$ is a featurizer
(assumed to be differentiable in its second argument).
We note that this predictor is linear in the output of the featurizer.

To select the model parameters, we solve a least squares problem with data given by
\[
A(\omega) =
\begin{bmatrix}
	e^{\omega^\mathrm{data}_1} \phi(u_1,\omega^\mathrm{feat})^T \\
	\vdots \\
	e^{\omega^\mathrm{data}_N} \phi(u_N,\omega^\mathrm{feat})^T \\
	e^{\omega^\mathrm{reg}_1} R_1 \\
	\vdots \\
	e^{\omega^\mathrm{reg}_d} R_d
\end{bmatrix}
, \qquad B(\omega) =
\begin{bmatrix}
	e^{\omega^\mathrm{data}_1} y_1 \\
	\vdots \\
	e^{\omega^\mathrm{data}_N} y_N \\
	0 \\
	\vdots \\
	0
\end{bmatrix},
\]
where $\omega^\mathrm{data}\in \Omega^\mathrm{data} \subseteq \reals^N$ are 
data weighting hyper-parameters, $R_1,\ldots,R_d$ are regularization matrices with 
appropriate sizes, and $\omega^\mathrm{reg} \in \Omega^\mathrm{reg} \subseteq \reals^d$
are regularization hyper-parameters.

The overall hyper-parameter is denoted
\[
\omega = (\omega^\mathrm{feat},\omega^\mathrm{data},\omega^\mathrm{reg})
\in \Omega^\mathrm{feat} \times \Omega^\mathrm{data} \times \Omega^\mathrm{reg}.
\]
We will describe the roles of each hyper-parameter in more detail below;
but here we note that $\omega^\mathrm{data}$ scales the individual training data examples,
$\omega^\mathrm{feat}$ are hyper-parameters in our featurizer, and
$\omega^\mathrm{reg}$ are hyper-parameters that scale each of our regularizers.
We assume that the hyper-parameter regularization function is separable, meaning it has the form
\[
r(\omega) = r^\mathrm{feat}(\omega^\mathrm{feat}) +
r^\mathrm{data}(\omega^\mathrm{data})
+ r^\mathrm{reg}(\omega^\mathrm{reg}),
\]
leading to a separable proximal operator \cite{parikh2014proximal}.

\subsection{True objective function}
We also have \emph{validation data} composed of inputs
$u^\mathrm{val}_1,\ldots,u^\mathrm{val}_{N_\mathrm{val}} \in \mathcal U$ and outputs 
$y^\mathrm{val}_1,\ldots,y^\mathrm{val}_{N_\mathrm{val}} \in \reals^m$.
We form predictions
\[
\hat y^\mathrm{val}_i = \phi(u^\mathrm{val}_i)^T \theta^\mathrm{ls}(\omega),
\]
where the featurizer is fixed, or $\phi(u) = \phi(u,\omega^\mathrm{data})$.
The true objective $\psi$ in least squares data fitting corresponds to the average loss
of our predictions of the validation outputs, which has the form
\[
\psi(\theta) = \frac{1}{N_\mathrm{val}} \sum_{i=1}^{N_\mathrm{val}} l(\hat y^\mathrm{val}_i, y^\mathrm{val}_i),
\]
where $l:\reals^m \times \reals^m \to \reals$ is a penalty function (assumed to be differentiable in its first argument).

\paragraph{Regression and classification.}
The setting that we have described encompasses many problems in data fitting, including both regression and classification.
In regression, the output is a scalar, \ie, $y\in\reals$.
In multi-task regression, the output is a vector, \ie, $y\in\reals^m$.
In boolean classification, $y\in \{e_1,e_2\}$ ($e_i$ is the $i$th unit vector in $\reals^2$),
and the output represents a boolean class.
In multi-class classification, $y\in\{e_i \mid i=1,\ldots,m\}$ ($e_i$ is the $i$th
unit vector in $\reals^m$), and the output represents a class label.

\paragraph{Regression penalty function.}
The penalty function in regression (and multi-task regression) problems often has the form
\[
l(\hat y, y) = \pi (r),
\]
where $r=\hat y - y$ is the residual and $\pi:\reals^m \to \reals$ is a penalty function applied to the residual.
Some common forms for $\pi$ are listed below.
\begin{itemize}
	\item \emph{Square}. The square penalty is given by $\pi(r) = \|r\|_2^2$.
	\item \emph{Huber}. The Huber penalty is a robust penalty that has the form of the 
square penalty for small residuals, and the 2-norm penalty for large residuals:
		\[
		\pi(r) = \begin{cases}\|r\|_2^2 & \|r\|_2 \leq M \\ M(2\|r\|_2-M) & \|r\|_2 > M.\end{cases}
		\]
        We can consider $M$ as a hyper-hyper-parameter.
	\item \emph{Bisquare}. The bisquare penalty is a robust penalty that is constant for large residuals:
		\[
		\pi(r) = \begin{cases}
		\frac{M^2}{6}\left(1-[1-(\frac{\|r\|_2^2}{M^2})]^3\right) & \|r\|_2 \leq M \\
		M^2/6 & \|r\|_2 > M.
		\end{cases}
		\]
		We can consider $M$ as a hyper-hyper-parameter.
\end{itemize}

\paragraph{Classification penalty function.}
For classification, we associate with our 
prediction $\hat y \in \reals^m$ the probability distribution
on the $m$ label values given by
\[
\Prob(y = e_i) = \frac{e^{\hat y_i}}{\sum_{j=1}^m e^{\hat y_j}}, \quad i=1,\ldots,m.
\]
We interpret our prediction $\hat y$ as giving us a distribution
on the labels of $y$, given $x$.  

We will use the \emph{cross-entropy loss} as the penalty function in 
classification.  It has the form
\[
l(\hat y, y = e_i) = -\hat y_i + \log(\sum_{j=1}^m e^{\hat y_j}), \quad i=1,\ldots,m.
\]
The true loss is then average negative log probability of $y$ under the predicted 
distribution, over the test set.

We now describe the role of each of the three (vector) components of our hyper-parameter
vector.

\subsection{Data weighting}
\label{sec:data-weighting}

We first describe the role of the data weighting
hyper-parameter $\omega^\mathrm{data}$.
The $i$th entry $\omega^\mathrm{data}_i$ \emph{weights} the squared error of the $(u_i,y_i)$
data point in the loss by $e^{2\omega^\mathrm{data}_i}$ (a positive number).
If $\omega_i$ is small, then the $i$th data point has little effect on the model parameter, and vice versa.

By separately weighting the loss values of each data point, we can get
the same effect as using a non-quadratic loss function.
However, instead of having to decide on which loss function to use, we
can automatically select the weights in a weighted square loss.

We let $\Omega^\mathrm{data} = \{x \mid \ones^T x = 0\}$,
meaning we constrain the geometric mean of $\exp(\omega^\mathrm{data})$ to be one.
We describe some forms for the hyper-parameter regularization function $r^\mathrm{data}$.
We can regularize the hyper-parameter towards each data point being weighted equally
($\omega^\mathrm{data}=0$), \eg, by using
$r^\mathrm{data}(\omega) = \lambda \|\omega\|_2^2$ or 
$r^\mathrm{data}(\omega) = \lambda\|\omega\|_1$,
where $\lambda>0$.  (Here $\lambda$ is a hyper-hyper-parameter, since
it scales a regularizer on the hyper-parameters.)

\paragraph{Proximal operator.}

We give details on a particular proximal operator that is needed later in the paper.
Evaluating the proximal operator of $r^\mathrm{data}(\omega) = \lambda \|\omega\|_2^2$
with $\Omega=\Omega^\mathrm{data}$ at $\nu$ with step size $t$
corresponds to solving the optimization problem
\BEQ
\begin{array}{ll}
\mbox{minimize} & t\lambda\|\omega\|_2^2 + (1/2)\|\omega-\nu\|_2^2 \\
\mbox{subject to} & \ones^T \omega = 0,
\end{array}
\nonumber
\EEQ
with variable $\omega$.
The (linear) KKT system for this optimization problem is
\[
\begin{bmatrix}
(1+\lambda t)I & \ones \\
\ones^T & 0
\end{bmatrix}
\begin{bmatrix}
\omega \\ y
\end{bmatrix}
=
\begin{bmatrix}
\nu \\
0
\end{bmatrix},
\]
with dual variable $y$.
This linear system can be solved efficiently using block elimination
\cite[Appendix C]{boyd2004convex}.

\subsection{Regularization}
\label{sec:reg}

The next hyper-parameter subvector that we consider is the regularization hyper-parameter $\omega^\mathrm{reg}$.
The regularization hyper-parameter affects the
\[
\sum_{i=1}^d \exp(2\omega^\mathrm{reg}_i) \|R_i \theta \|_F^2
\]
term in the least squares objective.
Each $\|R_i \theta \|_F^2$ term is meant to correspond to a measure of the complexity of $\theta$.
For example, if $R_i=I$, then the $i$th term is the
sum of the squares of the singular values of $\theta$.
The entries of the regularization hyper-parameter
correspond to the log of the weight on each regularization term.
The regularization matrices can have many forms; here we give two examples.

\paragraph{Diagonal regularization.}
Separate diagonal regularization has the form
\[
R_i = \diag(e_i), \quad i=1,\ldots,n,
\]
where $e_i$ is the $i$th unit vector in $\reals^{n}$.
The $i$th regularization term corresponds to the sum of squares of the $i$th row of $\theta$.

\paragraph{Graph regularization.}

The $R_i$ can correspond to incidence matrices of graphs between the elements in
each column of $\theta$.
Here the regularization hyper-parameter determines the relative
importance of the regularization graphs.

\subsection{Feature engineering}

The final hyper-parameter is the feature engineering hyper-parameter $\omega^\mathrm{feat}$, which parametrizes
the featurizer.
The goal is to select a $\omega^\mathrm{feat}$ which makes the output $y$ roughly
linear in $\phi(u,\omega^\mathrm{feat})$.
We assume that the input set $\mathcal U$ is a vector space in these examples.

\paragraph{Composition.}
Often $\phi$ is constructed as a feature generation chain, meaning it can be expressed as the composition
of individual feature engineering functions $\phi_1,\ldots,\phi_l$, or
\[
\phi = \phi_l \circ \cdots \circ \phi_1.
\]
Often the last feature engineering function adds a constant, or $\phi_l(x) = (x,1)$, so that the resulting predictor is affine.

\subsubsection{Scalar feature engineering functions}

We describe some scalar feature engineering functions $\phi:\reals \rightarrow \reals$, with the assumption that we could apply them
elementwise (with different hyper-parameters) to vector inputs.

\paragraph{Scaling.}
One of the simplest feature engineering functions is affine scaling, given by 
\[
\phi(x,(a,b)) = ax+b.
\]
It is common practice in data fitting to standardize or whiten the data, by scaling each
dimension with $a=1/\mathbf{std}(x)$ and $b=-\Expect[x]/\mathbf{std}(x)$.
Instead, with least squares auto tuning, we can select $a$ and $b$ based on the data.

\paragraph{Power transform.}

\begin{figure}
    \centering
    \includegraphics[width=.8\textwidth]{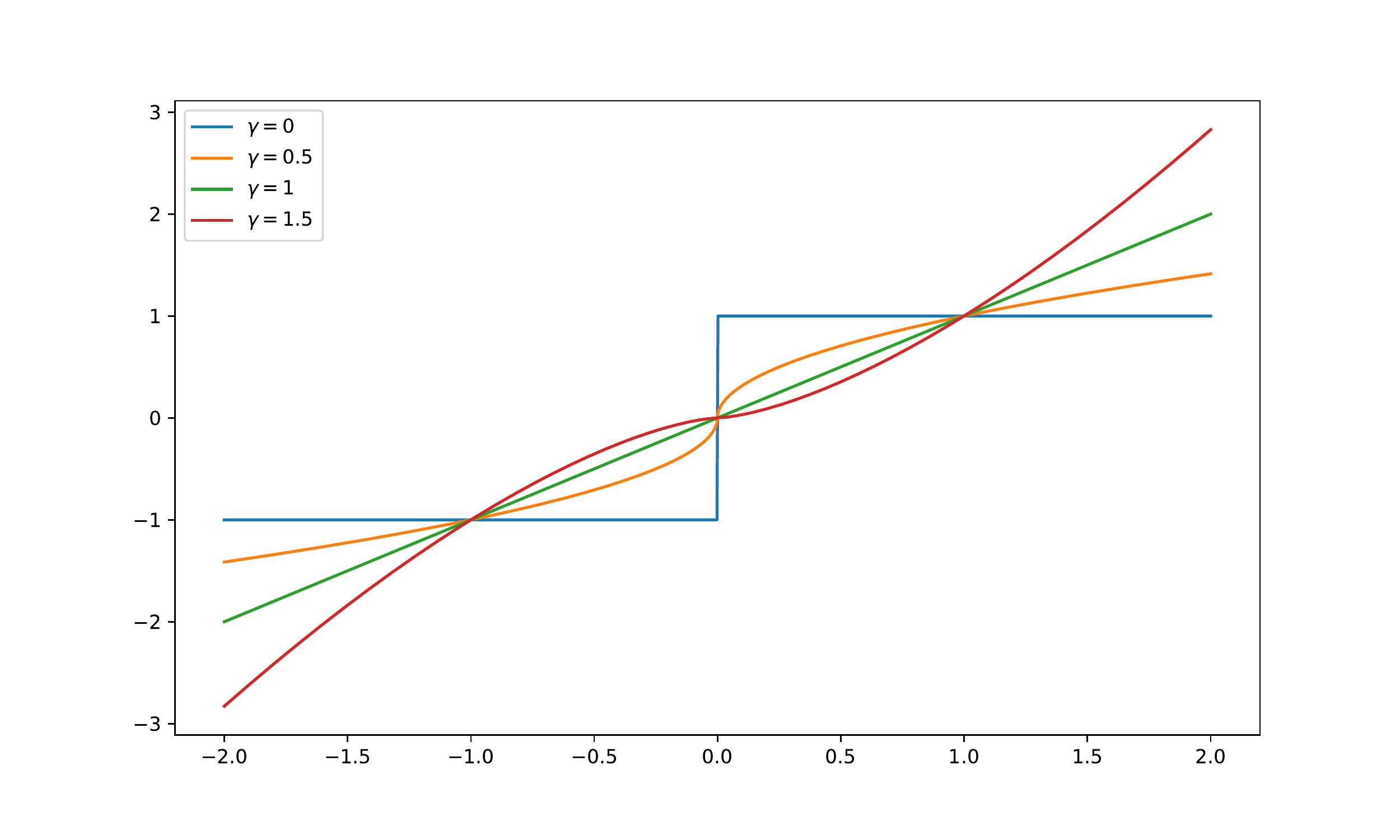}
    \caption{Examples of the power transform. Here the center $c=0$.}
    \label{fig:powers}
\end{figure}

The \emph{power transform} is given by
\[
\phi(x,(c,\gamma)) = {\bf sgn}(x-c)|x-c|^\gamma,
\]
where ${\bf sgn}(x)$ is $1$ if $x>0$, $-1$ if $x<0$ and $0$ if $x=0$, the center $c\in\reals$,
and the scale $\gamma \in\reals$.
Here the hyper-parameters are $\gamma\in\reals$ and the center $c\in\reals$.
For various values of $\gamma$ and $c$, this function defines
different transformations.
For example, if $\gamma=1$ and $c=0$, this transform is the identity.
If $\gamma=0$, this transform determines whether $x$ is to the right or
left of the center, $c$.
If $\gamma=1/2$ and $c=0$, this transform performs a symmetric square root.
This transform is differentiable everywhere except when $\gamma=0$.
See figure \ref{fig:powers} for some examples.

\paragraph{Polynomial splines.}

A spline is a piecewise polynomial function.
Given a monotonically increasing knot vector $z\in\reals^{k+1}$, a degree $d$, and polynomial
coefficients $f_0,\ldots,f_{k+1}\in\reals^{d+1}$, a spline is given by
\[
\phi(x,(z,f_1,\ldots,f_{k+1})) = \begin{cases}
\sum_{i=0}^d (f_0)_i x^i & x \in (-\infty,z_1) \\
\sum_{i=0}^d (f_j)_i x^i & x \in [z_j,z_{j+1}), \quad j=1,\ldots,k, \\
\sum_{i=0}^d (f_{k+1})_i x^i & x \in [z_{k+1},+\infty).
\end{cases}
\]
Here $r^\mathrm{feat}$ and $\Omega^\mathrm{feat}$ can be used to enforce continuity
(and differentiability) at $z_1,\ldots,z_{k+1}$.

\subsubsection{Multi-dimensional feature engineering functions}
Next we describe some multidimensional feature engineering functions.

\paragraph{Low rank.}
Can be a low rank transformation, given by
\[
\phi(x,T) = Tx
\]
where $T\in\reals^{r \times n}$, and $r < n$.
In practice, a common choice for $T$ is the first few eigenvectors of the singular value
decomposition of the data matrix.
With least squares auto tuning, we can select $T$ directly.

\paragraph{Neural networks.}

The featurizer $\phi$ can be a neural network;
in this case $\omega^\mathrm{feat}$ corresponds to the
neural network's parameters.

\paragraph{Feature selection.}
We can select a fraction $f$ of the features with
\[
\phi(x) = \diag(a) x,
\]
where $\Omega^\mathrm{feat} = \{a \mid a \in \{0,1\}^{n_1}, \ones^T a = \lfloor fn_1\rfloor\}$.
Here $f$ is a hyper-hyper-parameter.

\subsection{Test set and early stopping}
\label{sec:test}

There is a risk of overfitting to the validation set
when there are a large number of hyper-parameters, since we are
(almost) directly minimizing the validation loss.
To detect this, we introduce a third dataset, the \emph{test dataset},
and evaluate the fitted model on it \emph{once} at the end of the algorithm.
In particular, the validation loss throughout the algorithm need not be
an accurate measure of the model's performance on new unseen data,
especially when there are many hyper-parameters.

\paragraph{Early stopping.}
As a slight variation, to combat overfitting, we can
calculate the loss of the fitted model on the test at each iteration,
and halt the algorithm when the test loss begins to increase.
This technique is sometimes referred to as early stopping \cite{prechelt1998early}.
When performing early stopping, it is important to have a fourth dataset,
the \emph{final test dataset}, and evaluate on this set one time
when the algorithm terminates.
We have observed that this technique works very well in practice.
However, we do not use early stopping in our numerical example,
and instead run the algorithm until convergence.

\section{Numerical example}
\label{sec:num_ex}

In this section we apply our method of automatic least squares data fitting to
the well-studied MNIST handwritten digit classification dataset
\cite{lecun1998gradient}.
We note that in the machine learning community, this task is
considered ``solved'' by, \eg, deep convolutional neural networks (\ie, one can
achieve arbitrarily low test error).
We apply the ideas described in this paper to a large and small version of
MNIST in order to show that
standard least squares coupled with automatic tuning of additional
hyper-parameters can achieve relatively high test accuracy, and
can drastically improve the performance of standard least squares.
In this example, it is also worth nothing that we do not perform any
hyper-hyper-parameter optimization.

\paragraph{MNIST.}
The MNIST dataset is composed of 50,000 training data points,
where each data point is a 784-vector (a $28\times 28$ grayscale image
flattened in row-order).
There are $m=10$ classes, corresponding to the digits 0--9.
MNIST also comes with a test set, composed of 10,000 training points and labels.
Since the task here is classification, we use the cross-entropy loss as the true
objective function.
The code used to produce these results has been made freely available at
\url{www.github.com/sbarratt/lsat}.
All experiments were performed on an unloaded Nvidia 1080 TI GPU using floats.

We create two MNIST datasets by randomly selecting data points.
The \emph{small dataset} has 3,500 training data points and 1,500 validation data points.
The \emph{full dataset} has 35,000 training data points and 15,000 validation data points.
We evaluate four methods by tuning their hyper-parameters and then
calculating the final validation loss and test error.
The results are summarized in table \ref{tab:full} and table \ref{tab:small}.
We describe each method below, in order.

\begin{table}
  \caption{Small dataset.}
  \centering
  \begin{tabular}{llll}
    \toprule
    Method & Hyper-parameters & Validation loss & Test error (\%) \\
    \midrule
    LS 									   & 0    & 1.77 & 13.0 \\
    LS + reg $\times$ 2 				   & 2    & 1.76 & 11.6 \\
    LS + reg $\times$ 3 + feat 			   & 4    & 1.54 & 6.1  \\
    LS + reg $\times$ 3 + feat + weighting & 3504 & 1.54 & 6.0  \\
    \bottomrule
  \end{tabular}
  \label{tab:small}
\end{table}

\begin{table}
  \caption{Full dataset.}
  \centering
  \begin{tabular}{llll}
    \toprule
    Method & Hyper-parameters & Validation loss & Test error (\%) \\
    \midrule
    LS 									   & 0     & 1.74 & 10.3 \\
    LS + reg $\times$ 2 				   & 2     & 1.74 & 10.3 \\
    LS + reg $\times$ 3 + feat 			   & 4     & 1.53 & 4.7  \\
    LS + reg $\times$ 3 + feat + weighting & 35004 & 1.53 & 4.8  \\
    \bottomrule
  \end{tabular}
  \label{tab:full}
\end{table}

\paragraph{Base model.}

The simplest model is standard least squares, using the
$n=784$ image pixels as the feature vector.
That is, we solve the optimization problem
\BEQ
\begin{array}{ll}
\mbox{minimize} & \|X\theta - Y\|_F^2 + \|\theta\|_F^2.
\end{array}
\nonumber
\EEQ
Here we do no hyper-parameter tuning.
We refer to this model as \emph{LS} in the tables.

\paragraph{Regularization.}

To this simple model, we add a graph regularization term,
and optimize the two regularization hyper-parameters.
We define a graph on the length 784 feature vector, connecting
two nodes if the pixels they correspond to are adjacent to each
other in the original image.
We then compute the incidence matrix of this graph, as described in \S\ref{sec:reg},
and denote it by $A\in\reals^{1512 \times 784}$ (it has 1512 edges).
We then use the two regularization matrices:
\[
R_1 = I, \qquad R_2 = A.
\]
The matrix $R_1$ corresponds to standard ridge regularization, and $R_2$
measures the smoothness of the feature vector according to the graph we defined.
This introduces $2$ hyper-parameters, which separately weight $R_1$ and $R_2$.
We do not use a hyper-parameter regularization function, and initialize $\omega^\mathrm{reg}=(-2,-2)$.
We optimize these two regularization hyper-parameters to minimize validation loss.
We refer to this model as \emph{LS + reg $\times$ 2} in the tables.

\paragraph{Feature engineering.}

For each label, we run the $k$-means algorithm with $k=5$ on
the training data points that have that label.
From this, we get $km=50$ centers, which we call
\emph{archetypes} and denote by $a_1,\ldots,a_{50} \in \reals^{784}$.
Define the function $d$ such that it calculates how far $x$ is from
each of the archetypes, or
\[
d(x)_i = \|x-a_i\|_2, \quad i=1,\ldots,50.
\]
We use the feature engineering function
\[
\phi(x) = (x,s(-d(x)/\exp(\sigma)),1),
\]
where $\sigma$ is a feature engineering hyper-parameter,
and $s$ is the softmax function,
which transforms a vector $z\in\reals^n$ to a vector
in the probability simplex,
defined as
\[
s(z)_i = \frac{e^{z_i}}{\sum_j e^{z_j}}.
\]

We introduce separate ridge regularization for the
pixel features and the $k$-means features,
\ie, and still use the graph regularization on the pixel features;
the regularization matrices are
\[
R_1 = \begin{bmatrix}I & 0 & 0\end{bmatrix},
\quad R_2 = \begin{bmatrix}0 & I & 0\end{bmatrix},
\quad R_3 = \begin{bmatrix}A & 0 & 0\end{bmatrix}.
\]
We do not use a hyper-parameter regularization function for $\omega^\mathrm{feat}$.
We initialize $\sigma$ to $3$ and $\omega^\mathrm{reg}$ to $(0,0,0)$,
and optimize these four hyper-parameters.
We refer to this model as \emph{LS + reg $\times$ 3 + feat}.

\paragraph{Data weighting.}

\begin{figure}
    \centering
    \begin{subfigure}[b]{0.16\textwidth}
        \includegraphics[width=\textwidth]{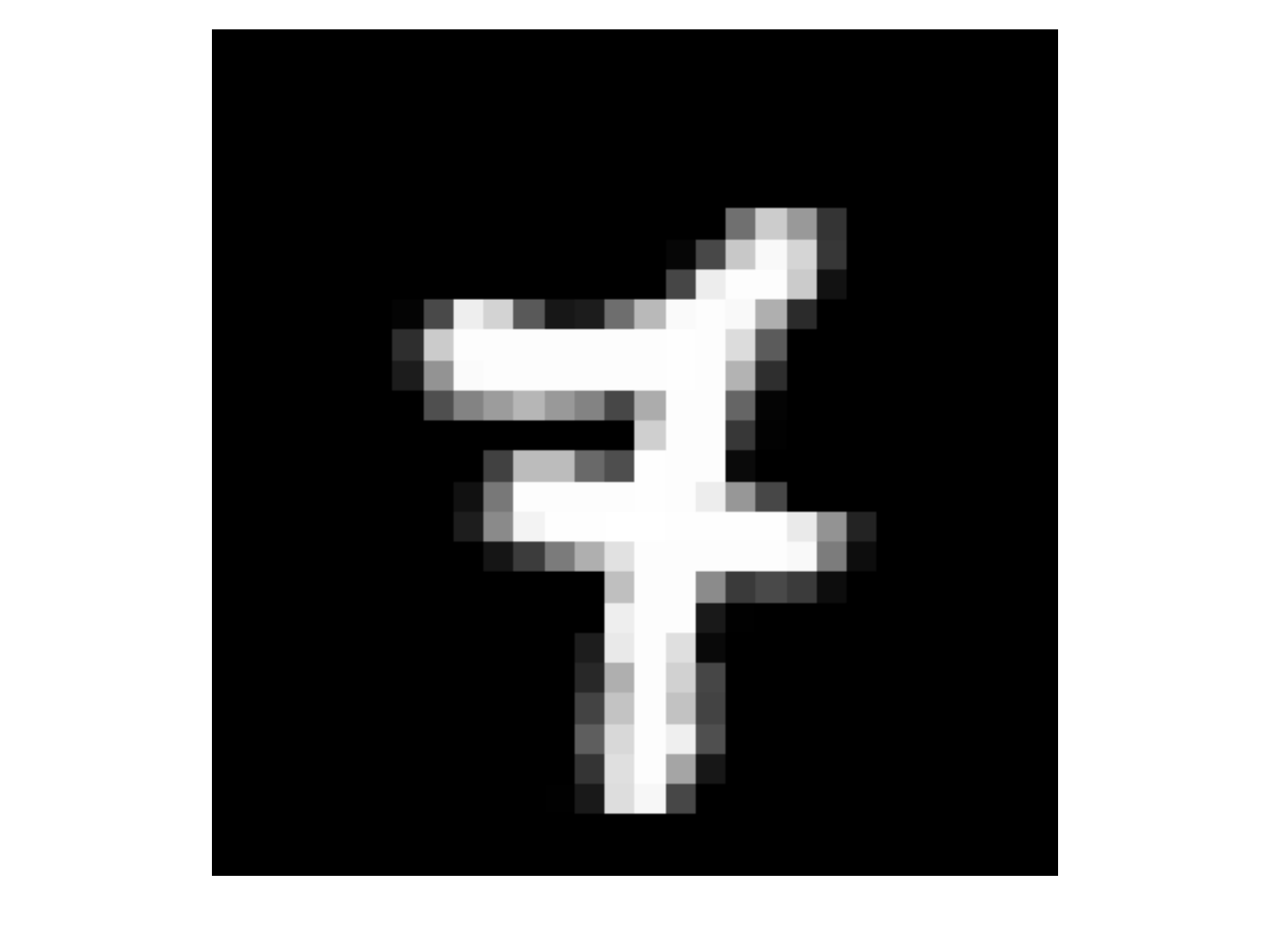}
        \caption{7}
    \end{subfigure}
    \hfill
    \begin{subfigure}[b]{0.16\textwidth}
        \includegraphics[width=\textwidth]{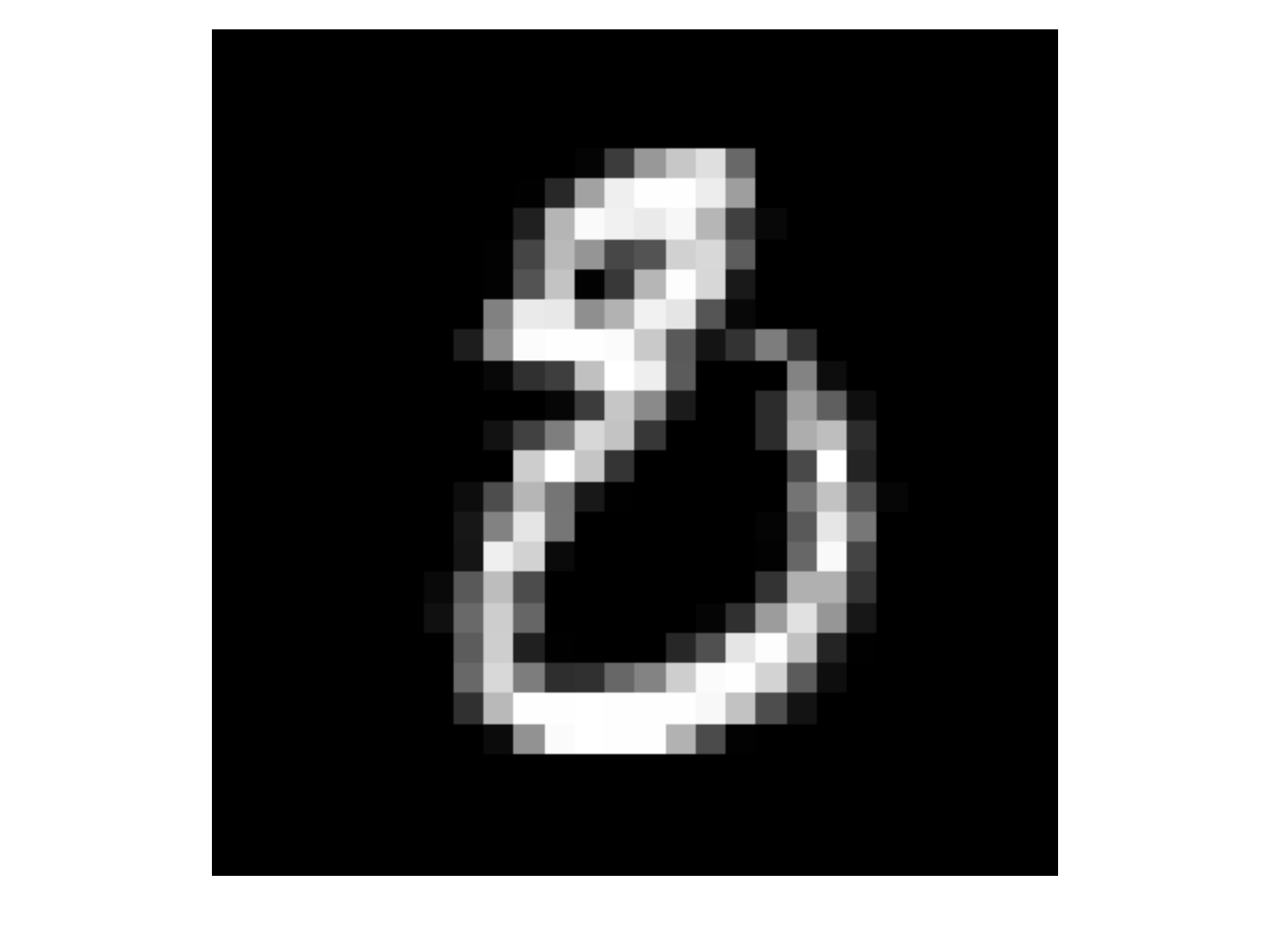}
        \caption{8}
    \end{subfigure}
    \hfill
    \begin{subfigure}[b]{0.16\textwidth}
        \includegraphics[width=\textwidth]{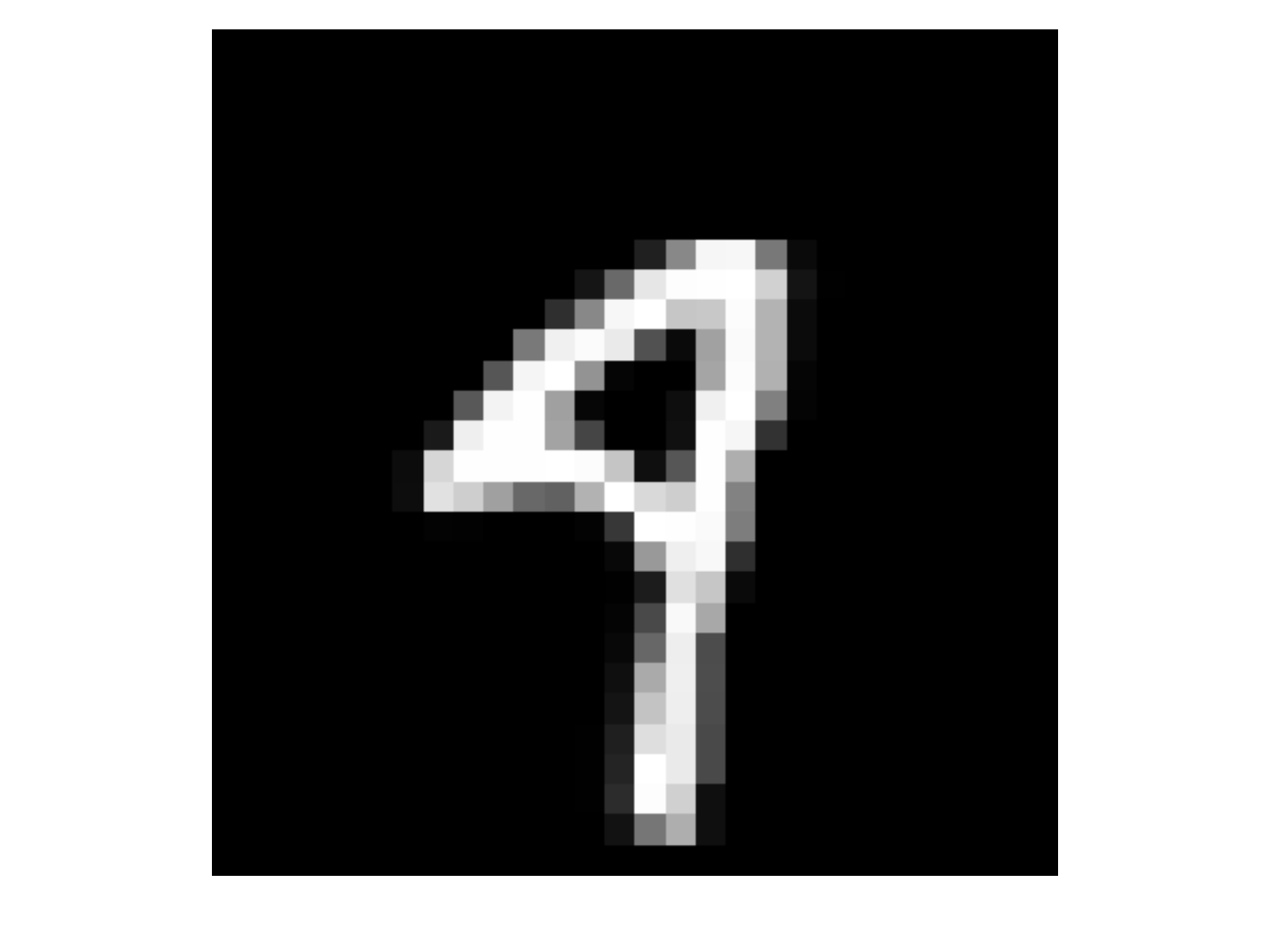}
        \caption{4}
    \end{subfigure}
    \hfill
    \begin{subfigure}[b]{0.16\textwidth}
        \includegraphics[width=\textwidth]{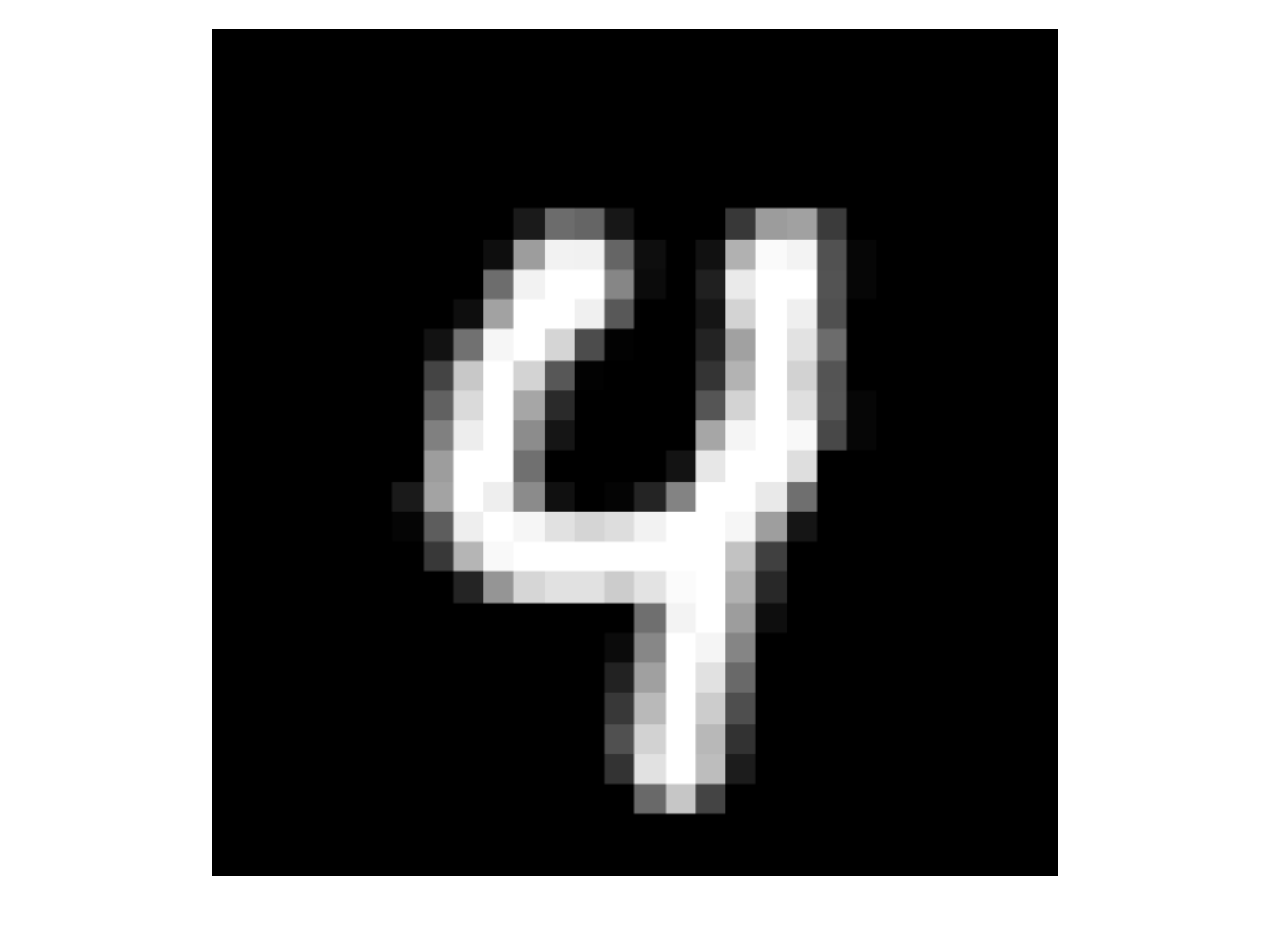}
        \caption{4}
    \end{subfigure}
    \hfill
    \begin{subfigure}[b]{0.16\textwidth}
        \includegraphics[width=\textwidth]{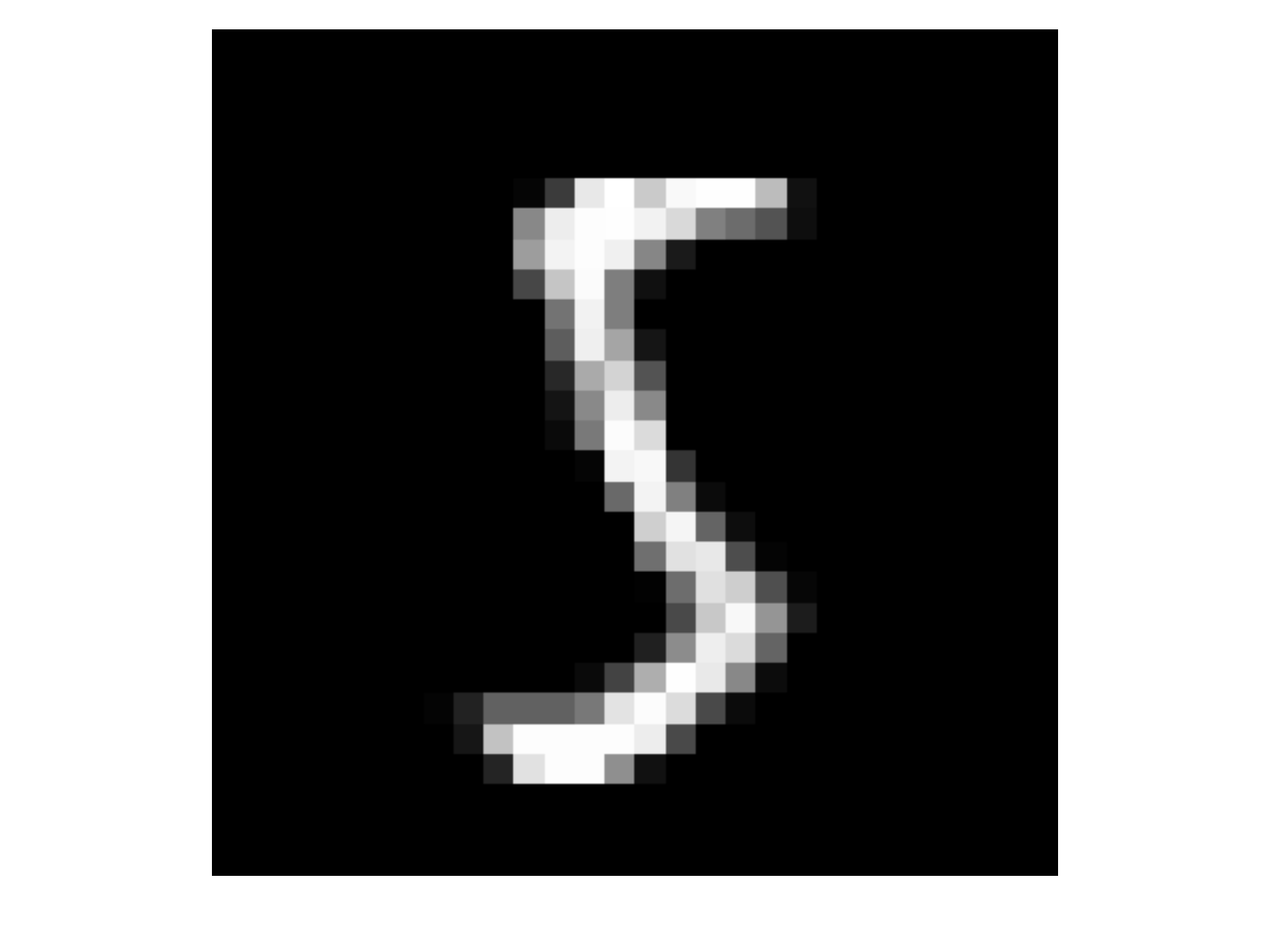}
        \caption{5}
    \end{subfigure}
    \hfill
    \begin{subfigure}[b]{0.16\textwidth}
        \includegraphics[width=\textwidth]{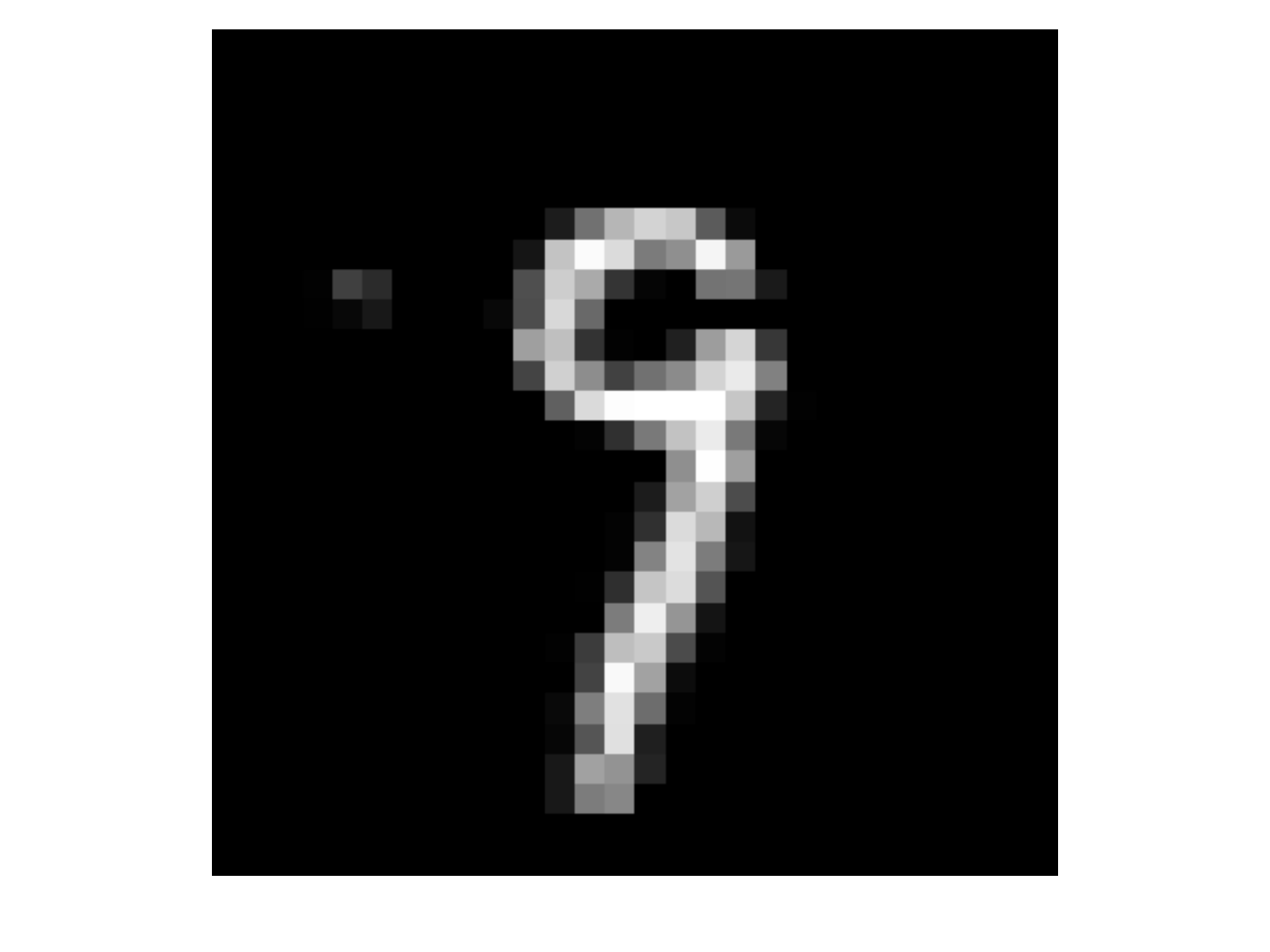}
        \caption{9}
    \end{subfigure}
    \caption{Training data with the lowest weights. Captions correspond to labels.}
    \label{fig:outliers}
\end{figure}

\begin{figure}
    \centering
    \begin{subfigure}[b]{0.16\textwidth}
        \includegraphics[width=\textwidth]{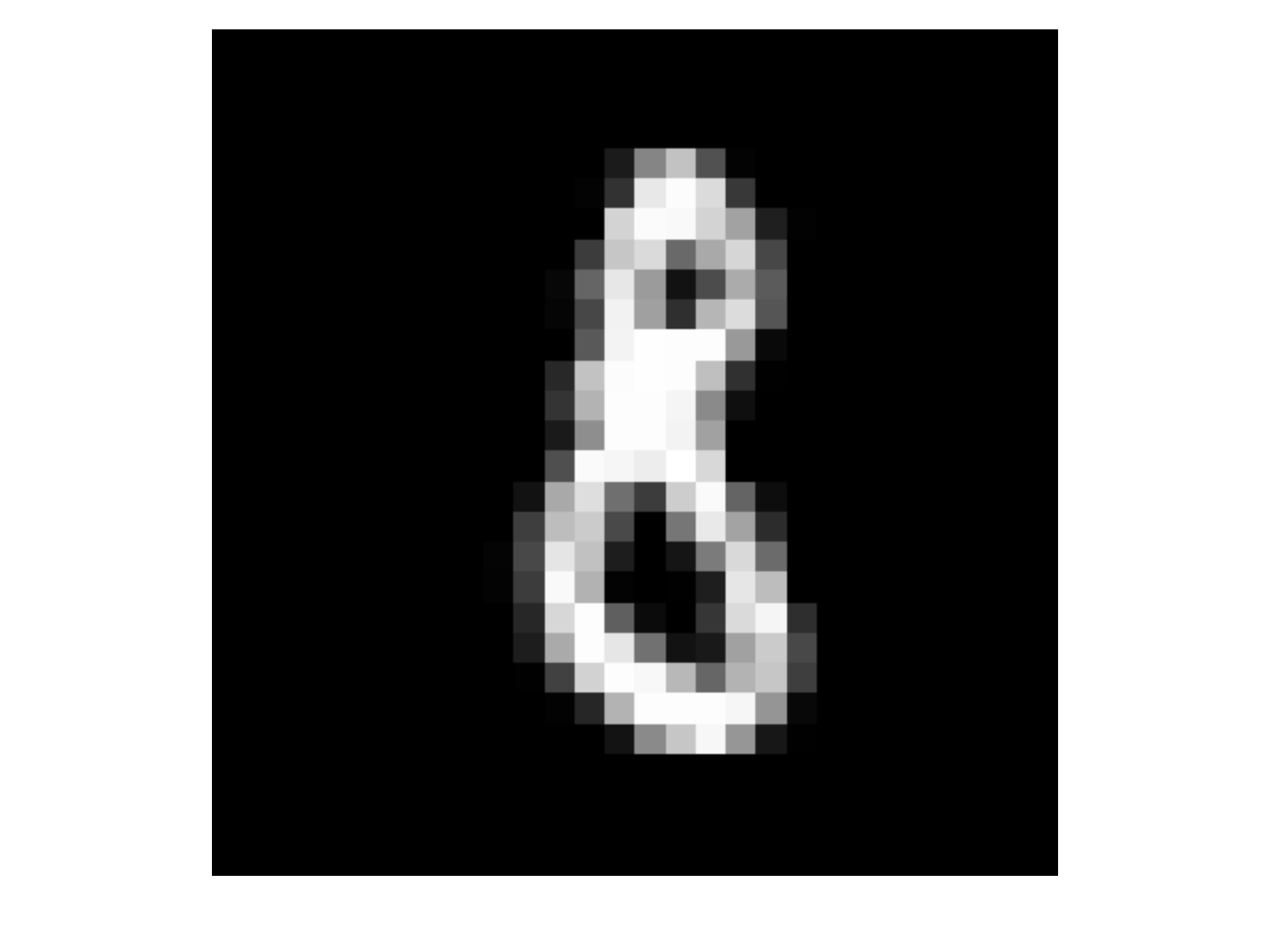}
        \caption{8}
    \end{subfigure}
    \hfill
    \begin{subfigure}[b]{0.16\textwidth}
        \includegraphics[width=\textwidth]{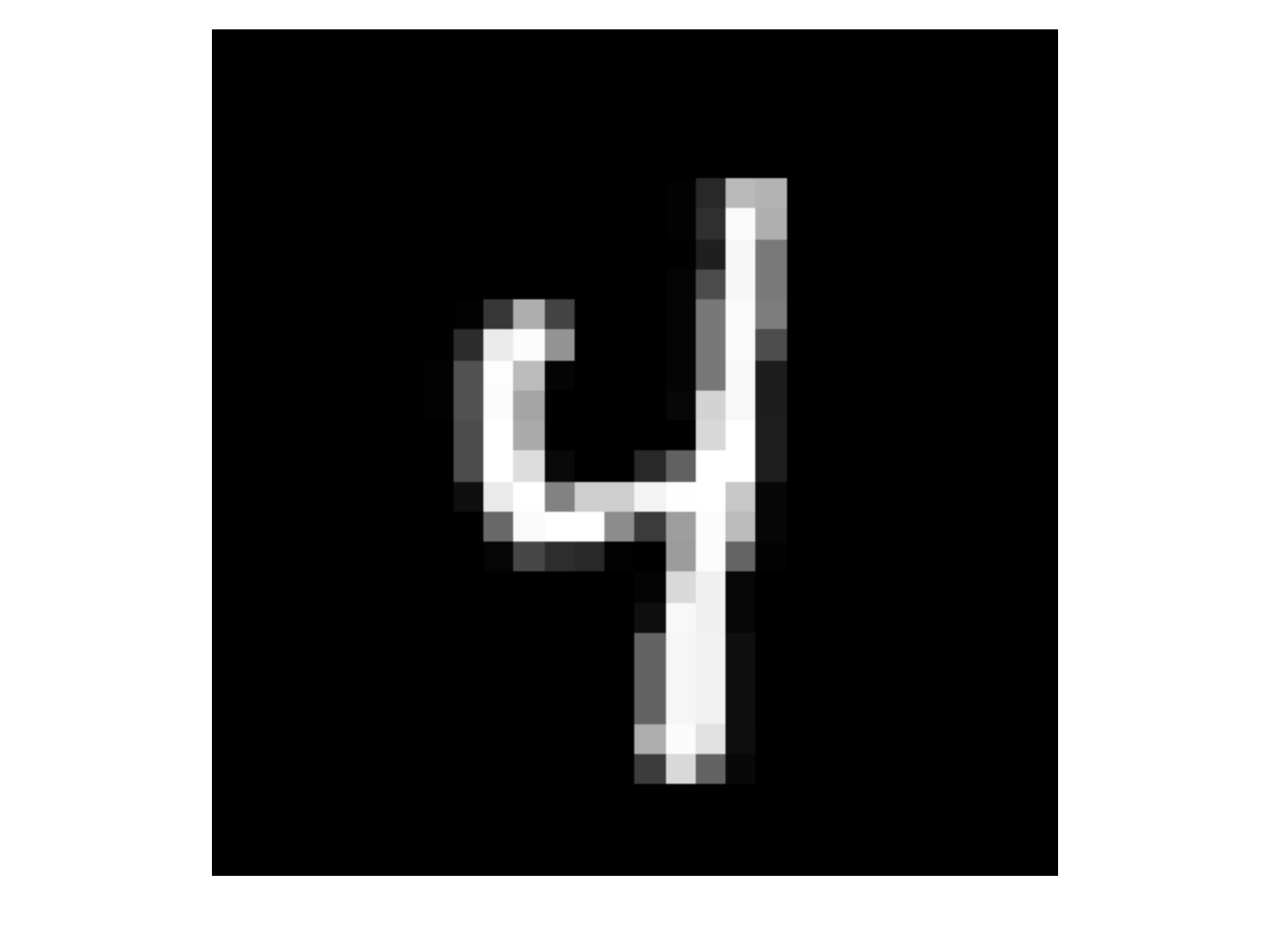}
        \caption{4}
    \end{subfigure}
    \hfill
    \begin{subfigure}[b]{0.16\textwidth}
        \includegraphics[width=\textwidth]{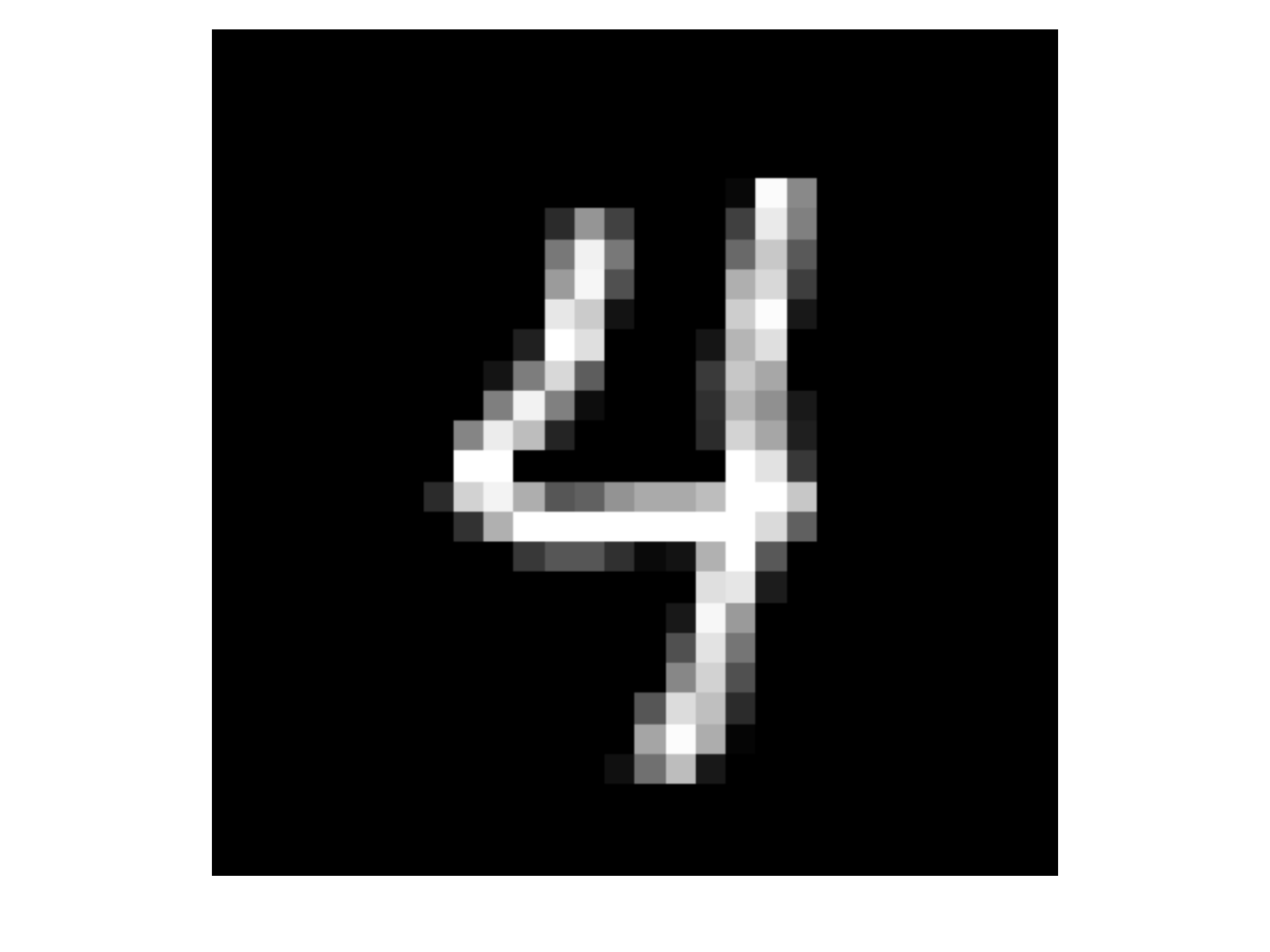}
        \caption{4}
    \end{subfigure}
    \hfill
    \begin{subfigure}[b]{0.16\textwidth}
        \includegraphics[width=\textwidth]{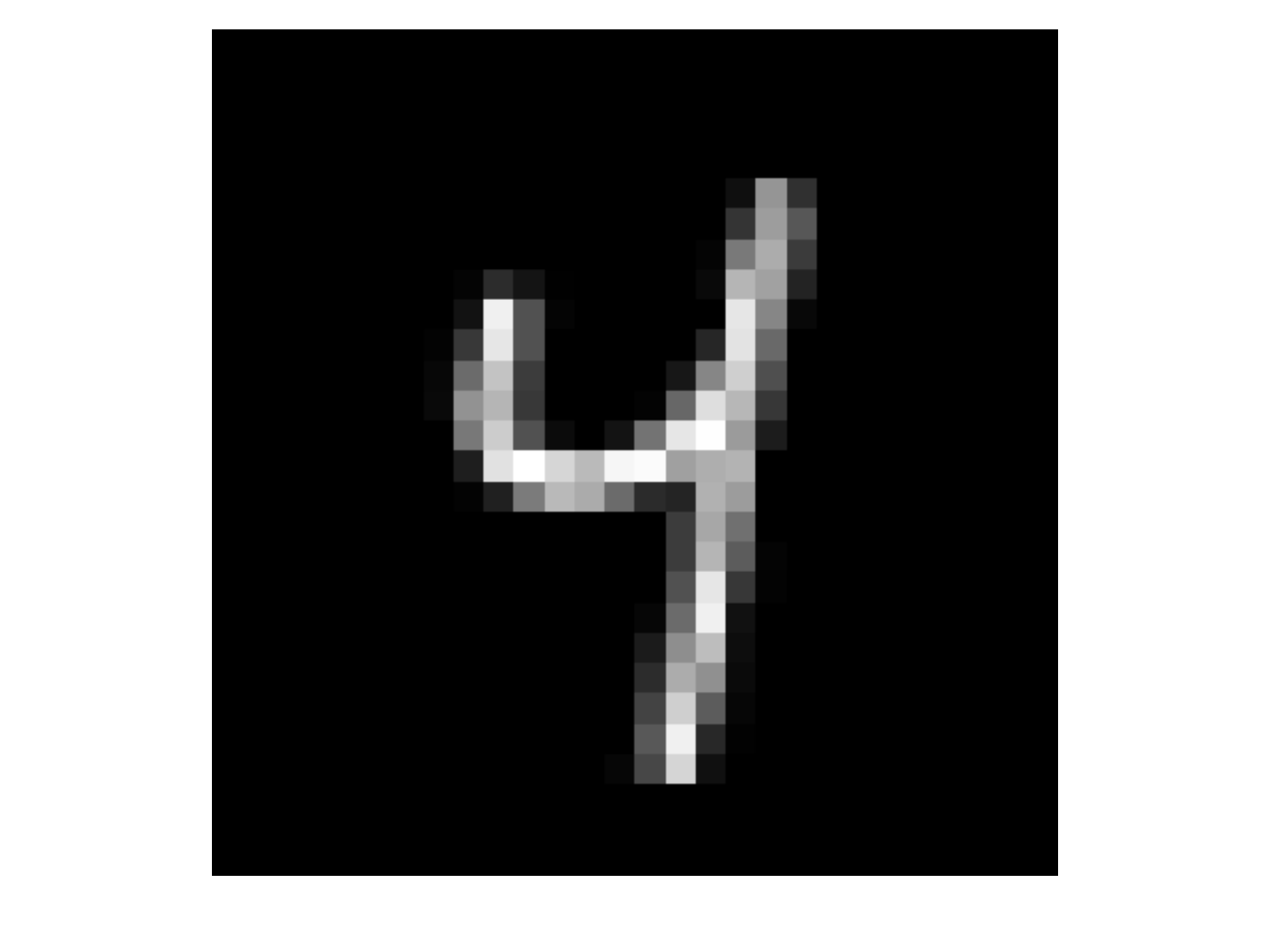}
        \caption{4}
    \end{subfigure}
    \hfill
    \begin{subfigure}[b]{0.16\textwidth}
        \includegraphics[width=\textwidth]{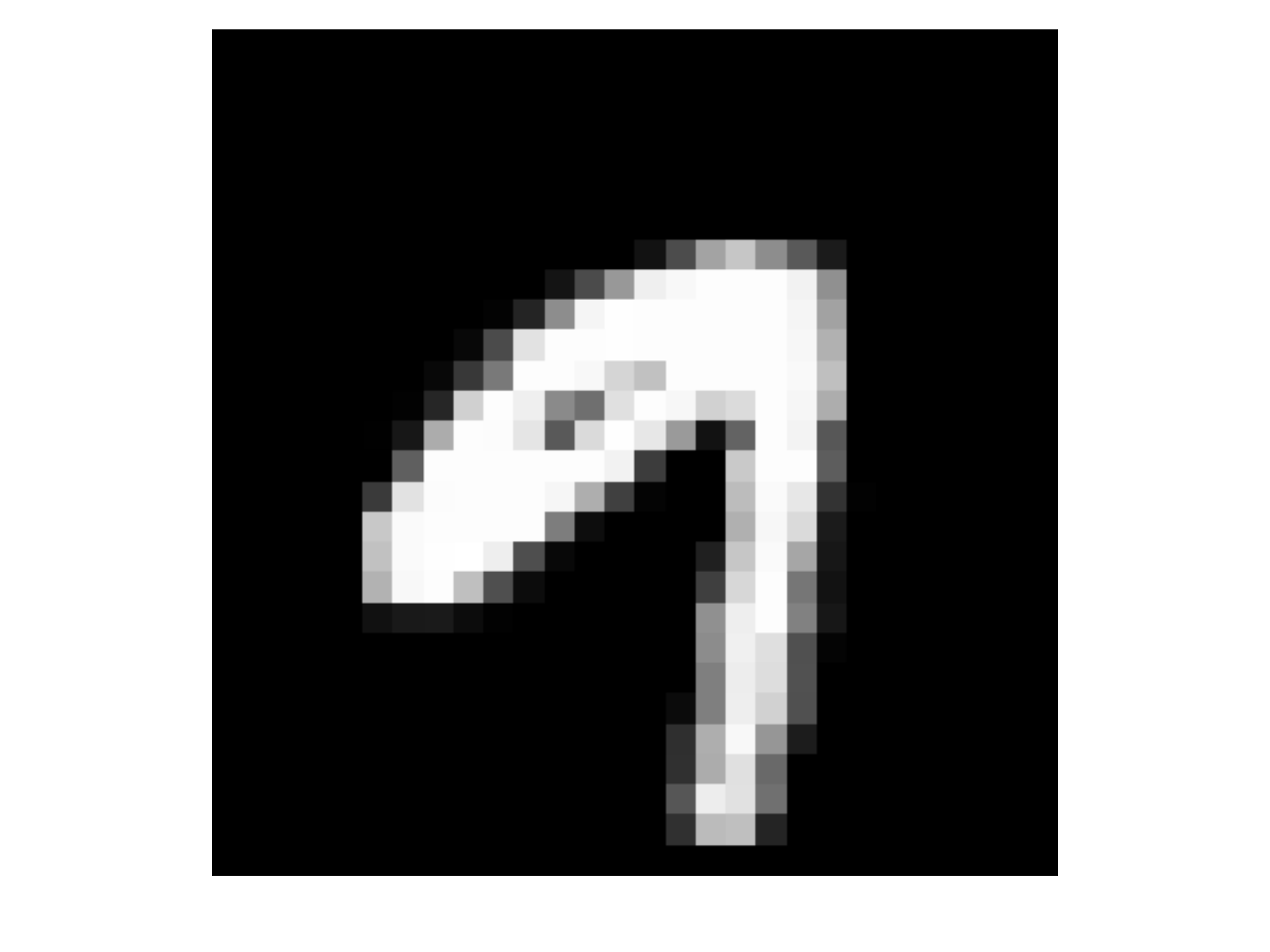}
        \caption{9}
    \end{subfigure}
    \hfill
    \begin{subfigure}[b]{0.16\textwidth}
        \includegraphics[width=\textwidth]{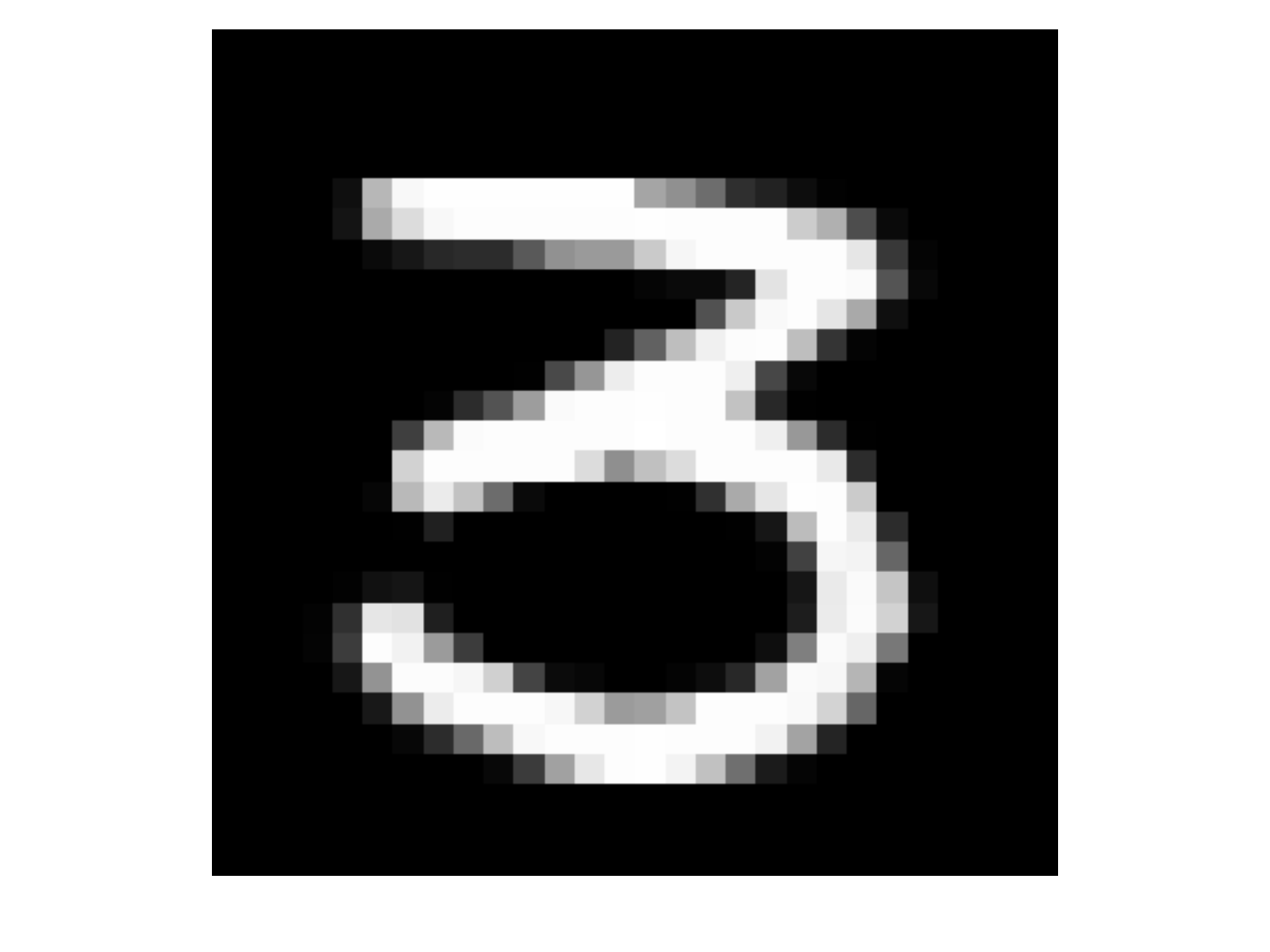}
        \caption{3}
    \end{subfigure}
    \caption{Training data with the highest weights. Captions correspond to labels.}
    \label{fig:exemplars}
\end{figure}

To LS + reg $\times$ 3 + feat, we add data weighting,
as described in \S\ref{sec:data-weighting}.
We use $\Omega^\mathrm{reg} = \{\omega \mid \ones^T \omega = 0\}$
and $r^\mathrm{reg}(\omega^\mathrm{reg}) = (0.01)\|\omega^\mathrm{reg}\|_2^2$.
This introduces 3,500 hyper-parameters in the case of the small dataset,
and 35,000 hyper-parameters for the large dataset.
We use the initialization $\omega=0$.
We refer to this model as \emph{LS + reg $\times$ 3 + feat + weighting}.
This method performs the best on the small dataset, in terms of test error.
On the full dataset, it performs slightly worse than the model without data
weighting, likely because of the overfitting phenomenon
discussed in \S\ref{sec:test}.
We show the training examples with the lowest data weights
and the training examples with the highest weights in
figure \ref{fig:outliers} and figure \ref{fig:exemplars} respectively,
on the small dataset.
The training data points with low weights seem harder to classify
(for example, (b) and (c) in figure \ref{fig:outliers} could be interpreted as
nines).

\section{Conclusion}

The authors are currently writing a second paper,
\emph{Least Squares Auto-Tuning Examples}, which
will detail many more applications of the methods
described in this paper to data fitting, control, and estimation.

\clearpage
\section*{Acknowledgements}
Shane Barratt is supported by the National Science Foundation Graduate Research Fellowship
under Grant No. DGE-1656518.

\bibliography{lsat}

\appendix

\clearpage
\section{Derivation of gradient of least squares solution}
\label{sec:deriv_grads}

Consider the least squares solution map
$\phi: \reals^{k \times n} \times \reals^{k \times m} \to \reals^{n \times m}$, given by
\[
\theta=\phi(A,B)=(A^TA)^{-1}A^TB.
\]
We are interested in the linear operator $D_A \phi(A,B): \reals^{k \times n} \rightarrow \reals^{n \times m}$,
\ie, the derivative of $\phi$ with respect to $A$, and the linear operator
$D_B \phi(A,B): \reals^{k \times m} \rightarrow \reals^{n \times m}$,
\ie, the derivative of $\phi$ with respect to $B$.

\paragraph{Derivative with respect to $A$.}
We have that
\[
D_A\phi(A,B)(\Delta A) = (A^TA)^{-1}\Delta A^T B - (A^TA)^{-1}(\Delta A^T A+A^T\Delta A) \theta,
\]
since
\BEASn
\phi(A+\Delta A,B) &=& ((A+\Delta A)^T(A+\Delta A))^{-1}(A^TB + \Delta A^TB) \\
&\approx& (A^TA)^{-1}\left(I-\Delta A^T A(A^TA)^{-1}-A^T\Delta A(A^TA)^{-1}\right)(A^TB + \Delta A^TB) \\
&\approx& \phi(A,B) + (A^TA)^{-1}\Delta A^T B - (A^TA)^{-1}(\Delta A^T A+A^T\Delta A) \theta,
\EEASn
where we used the approximation $(X+Y)^{-1} \approx X^{-1}-X^{-1}YX^{-1}$ for $Y$ small, and dropped higher order terms.
Suppose $f=\psi \circ \phi$ and $C=(A^TA)^{-1}\nabla_{\theta} \psi$
for some $\psi: \reals^{n \times m} \to \reals$.
Then the linear map $D_A f(A,B)$ is given by
\BEASn
D_A f (A,B) (\Delta A) &=& D_{\theta^{\mathrm{ls}}} \psi (D_A \phi(A,B) (\Delta A)) \\
&=& \Tr\left(\nabla_{\theta} \psi^T \left((A^TA)^{-1}\Delta A^T B - (A^TA)^{-1}(\Delta A^T A+A^T\Delta A) \theta\right)\right) \\
&=& \Tr((BC^T-A\theta C^T-AC\theta^T)^T\Delta A),
\EEASn
from which we conclude that $\nabla_A f = (B-A\theta)C^T-AC\theta^T$.

\paragraph{Derivative with respect to $B$.}
We have that $D_B\phi(A,B)(\Delta B) = (A^TA)^{-1}A^T\Delta B$, since
\BEASn
\phi(A,B+\Delta B) = \phi(A,B) + (A^TA)^{-1}A^T\Delta B.
\EEASn
Suppose $f=\psi \circ \phi$ for some $\psi: \reals^{n \times m} \to \reals$
and $C=(A^TA)^{-1}\nabla_{\theta} \psi$.
Then the linear map $D_B f(A,B)$ is given by
\BEASn
D_B f(A,B) (\Delta B) &=& D_{\theta^{\mathrm{ls}}} \psi (D_B \phi(A,B) (\Delta B)) \\
&=& \Tr(\nabla_{\theta} \psi^T (A^TA)^{-1}A^T\Delta B) \\
&=& \Tr((AC)^T\Delta B),
\EEASn
from which we conclude that $\nabla_B f = AC$.

\section{Derivation of stopping criterion}
\label{sec:deriv_stopping_criterion}

The optimality condition for minimizing $\psi + r$ is
\[
\nabla_\omega \psi + g = 0,
\]
where $g\in\partial_\omega r$, the subdifferential of $r$.
We have that
\[
t^k\partial r(\omega^{k+1}) + \omega^{k+1} - \omega^k + t^k \nabla_{\omega^k} \psi = 0,
\]
which implies that
\[
(\omega^k-\omega^{k+1})/t^k - g^k \in \partial r(\omega^{k+1}).
\]
Therefore, the optimality condition for $\omega^{k+1}$ is
\[
(\omega^k-\omega^{k+1})/t^k + (g^{k+1} - g^k) = 0,
\]
which leads to the stopping criterion \eqref{eq:stopping-criterion}.

\section{PyTorch implementation}

We implemented the forward/backward methods of the least squares in PyTorch.
We compute the forward pass using a Cholesky factorization of the Gram matrix $A^TA$,
which is cached and reused in the backward pass.
The code is simply:
\begin{verbatim}
import torch

class DenseLeastSquares(torch.autograd.Function):
    @staticmethod
    def forward(ctx, A, B):
        with torch.no_grad():
            u = torch.cholesky(A.t() @ A, upper=True)
            theta = torch.potrs(A.t() @ B, u)
        ctx.save_for_backward(A, B, theta, u)

        return theta

    @staticmethod
    def backward(ctx, dtheta):
        A, B, theta, u = ctx.saved_tensors

        with torch.no_grad():
            C = torch.potrs(dtheta, u)
            Cthetat = C @ theta.t()
            dA = B @ C.t() - A @ (Cthetat + Cthetat.t())
            dB = A @ C

        return dA, dB
\end{verbatim}

\section{Tensorflow implementation}

We implemented the forward/backward methods of least squares in Tensorflow.
We compute the forward pass using a Cholesky factorization of the Gram matrix $A^TA$, which
is cached and reused in the backward pass.
The code is simply:
\begin{verbatim}
import tensorflow as tf

@tf.custom_gradient
def lstsq(A,B):
    AtA = tf.transpose(A)@A
    chol = tf.linalg.cholesky(AtA)
    theta = tf.linalg.cholesky_solve(chol, tf.transpose(A)@B)
    def grad(dtheta):
        C = tf.linalg.cholesky_solve(chol, dtheta)
        Cthetat = C@tf.transpose(theta)
        dA = B@tf.transpose(C)-A@(Cthetat+tf.transpose(Cthetat))
        dB = A@C
        return dA, dB
    return theta, grad
\end{verbatim}

\end{document}

%% file: defs.tex
\newcommand{\BEAS}{\begin{eqnarray}}
\newcommand{\EEAS}{\end{eqnarray}}
\newcommand{\BEASn}{\begin{eqnarray*}}
\newcommand{\EEASn}{\end{eqnarray*}}
\newcommand{\BEQ}{\begin{equation}}
\newcommand{\EEQ}{\end{equation}}
\newcommand{\BIT}{\begin{itemize}}
\newcommand{\EIT}{\end{itemize}}

\newcommand{\eg}{{\it e.g.}}
\newcommand{\ie}{{\it i.e.}}

\newcommand{\ones}{\mathbf 1}
\newcommand{\reals}{{\mbox{\bf R}}}

\newcommand{\Tr}{\mathop{\bf tr}}

\newcommand{\diag}{\mathop{\bf diag}}

\newcommand{\Expect}{\mathop{\bf E{}}}
\newcommand{\Prob}{\mathop{\bf Prob}}

\newcommand{\argmin}{\mathop{\rm argmin}}

\newcounter{algorithmctr}[section]
\renewcommand{\thealgorithmctr}{\thesection.\arabic{algorithmctr}}
\newenvironment{algdesc}%
   {\refstepcounter{algorithmctr}\begin{list}{}{%
       \setlength{\rightmargin}{0\linewidth}%
       \setlength{\leftmargin}{.05\linewidth}}%
       \rmfamily\small
       \item[]{\setlength{\parskip}{0ex}\hrulefill\par%
        \nopagebreak{\bfseries\textsf{Algorithm \thealgorithmctr~}}}}%
   {{\setlength{\parskip}{-1ex}\nopagebreak\par\hrulefill} \end{list}}